\documentclass[11pt]{article}
\usepackage{hyperref}
\usepackage{graphicx,amsfonts,amssymb,amsmath,amsthm}
\usepackage{ucs}
\usepackage{array}
\usepackage{url}
\usepackage{colortbl}
\usepackage[utf8x]{inputenc}
\usepackage[T1]{fontenc}
\usepackage{ulem}

\usepackage{multirow}
\usepackage[hmargin=2cm,vmargin=2cm]{geometry}
\begin{document}\bibliographystyle{ieeetr}
\normalem

\title{Parsimonious Mahalanobis Kernel for the Classification of High Dimensional Data}%

\author{M.~Fauvel$^1$, A.~Villa$^2$, J.~Chanussot$^3$ and J.A.~Benediktsson$^4$}
\date{$^1$ INRA, DYNAFOR, BP 32607 - Auzeville-Tolosane 31326 - Castanet Tolosan - FRANCE\\
  $^2$ Aresys srl, via Bistolfi 49, 20134 Milano and Dipartimento di Elettronica ed Informazione, Politecnico di Milano, 20133 Milano, Italy\\
  $^3$ GIPSA-lab, Departement Image Signal, BP 46 - 38402 Saint Martin d'H\`eres - FRANCE\\
  $^4$ Dept. of Electrical and Computer Engineering, University of Iceland Hjardarhagi 2-6, 107 Reykjavik - ICELAND  
}

\maketitle

\begin{abstract}
  The classification  of high dimensional data with  kernel methods is
  considered  in   this  article.   Exploiting   the  \emph{emptiness}
  property  of  high  dimensional   spaces,  a  kernel  based  on  the
  Mahalanobis distance is proposed. The computation of the Mahalanobis
  distance  requires the inversion  of a  covariance matrix.   In high
  dimensional   spaces,    the   estimated   covariance    matrix   is
  ill-conditioned and its inversion  is unstable or impossible.  Using
  a parsimonious statistical  model, namely the \emph{High Dimensional
    Discriminant  Analysis}  model,  the  specific  signal  and  noise
  subspaces are estimated for each considered class making the inverse
  of the class specific covariance matrix explicit and stable, leading
  to the  definition of  a \emph{parsimonious Mahalanobis}  kernel.  A
  SVM based framework is used for selecting the hyperparameters of the
  parsimonious   Mahalanobis  kernel   by  optimizing   the  so-called
  \emph{radius-margin}  bound.   Experimental  results on  three  high
  dimensional data sets show that  the proposed kernel is suitable for
  classifying high  dimensional data, providing  better classification
  accuracies than the conventional Gaussian kernel.
\end{abstract}

\section{Introduction}
High Dimensional  (HD) data sets  are commonly available for  fully or
partially  automatic  processing:  For  a  relatively  low  number  of
samples,  $n$, a  huge  number of  variables,  $d$, is  simultaneously
accessible.   For  instance,  in  hyperspectral imagery,  hundreds  of
spectral  wavelengths  are  recorded   for  a  given  pixel;  in  gene
expression analysis,  the measure of expression level  of thousands of
genes is typical; in customer recommendation systems for web services,
to each potential client a high number of variables is associated (his
past         choices,         his         personal         information
...)~\cite{Kriegel:2009:CHD:1497577.1497578,donoho,Plaza20041097}.
For  each  sample,  it  is   possible  to  have  either  numerical  or
alphabetical variables which can be  sparse or with a different signal
to noise ratio.

In terms  of processing, such data  may need to  be either classified,
clustered, filtered  or inversed in a supervised  or unsupervised way.
Although many algorithms exist in the literature for small or moderate
dimensions (from Bayesian methods to Machine Learning techniques) most
of  them are  not well  suited  to HD  data.  Actually,  HD data  pose
critical theoretical and practical  problems that need to be addressed
specifically~\cite{donoho}.

\textcolor{black}{Indeed, HD  spaces exhibit non  intuitive geometrical
  and  statistical  properties  when  compared  to  lower  dimensional
  spaces. Most  of them  do not behave  in a similar  way as  in three
  dimensional Euclidean  spaces (Table~\ref{table:hd:space} summarizes
  the   main  properties   of  HD   spaces)~\cite{book:kendall}.   For
  instance, samples  following a uniform  law will have a  tendency to
  have a  high concentration in  the corners~\cite{jimenez:landgrebe}.
  The same property holds for  normally distributed data: samples tend
  to  have a  high  concentration  in the  tails~\cite{herault:esann},
  making density  estimation a  difficult task.   This problem  can be
  related to  the number of  parameters $t$ to  be estimated to  fit a
  Gaussian  distribution  which  grows quadratically  with  the  space
  dimensionality, $t=d(d+3)/2$  (5150 for $d=100$).  Because  of this,
  conventional generative methods are  not suitable for analyzing this
  type of data.}

\textcolor{black}{Unfortunately, discriminative methods  also suffer if
  the  dimensionality is  high,  due to  the ``\emph{concentration  of
    measure phenomenon}''~\cite{donoho}.   In HD spaces,  samples tend
  to        be         equally        distant         from        each
  other~\cite{springerlink:10.1007/3-540-44503-X_27}.   Hence,  it  is
  clear  that  Nearest Neighbors  methods  will  definitively fail  to
  process such  data.  Moreover,  the Euclidean  distance will  not be
  appropriate to assess the similarity  between two samples.  In fact,
  it    is    has   been    shown    that    every   Minkowski    norm
  ($\|\mathbf{x}\|_m=\big(\sum_{i=1}^d|\mathbf{x}_i|^m   \big)^{1/m}$,
  $m=1,        2\dotsc$)       is        affected       by        this
  phenomenon~\cite{Francois:2007:CFD:1271928.1272109}.      Therefore,
  every method  based on  the distance  between samples~\cite{1316859}
  (SVM  with  Gaussian  kernel,  neural  network,  Nearest  Neighbors,
  Locally  Linear Embedding\dots)  are  potentially  affected by  this
  phenomenon~\cite{springerlink:10.1007/3-540-44869-1_14,verleysen}.}

An additional property, for  which the consequences are more practical
than      theoretical,       is      the      ``\emph{empty      space
  phenomenon}''~\cite{citeulike:4303943}: In  HD spaces, the available
samples usually fill  a very small part of  the space. Therefore, most
of  the space  is  empty.  Note  that  if originally  the empty  space
phenomenon  was  considered as  a  problem, it  will  be  seen in  the
following that it is actually  the basis of several useful statistical
models.

Today,  the phrasing  ``\emph{curse  of dimensionality}'',  originally
from R.~Bellman~\cite{citeulike:4303943}, refers to the aforementioned
problems of HD data and  reflects how processing HD data is difficult.
However,  as  D.~Donoho has  noticed~\cite{donoho},  there  is also  a
``\emph{Blessing of dimensionality}'': For instance in classification,
the class separability is improved when the dimensionality of the data
increases.   Consider for example  a comparison  between hyperspectral
(hundreds of spectral wavelengths) and multispectral (tens of spectral
wavelengths)  remote sensing images\cite{book:landgrebe}.   The former
contains  much   more  information,   and  enables  a   more  accurate
distinction  of  the land  cover  classes.   However, if  conventional
methods   are   used,   the   additional  information   contained   in
hyperspectral   images  will   not  lead   to  an   increase   of  the
classification    accuracy~\cite{jimenez:landgrebe}.    Hence,   using
conventional methods, classification accuracies remains low.

Several methods have  been proposed in the literature  to deal with HD
data for  the purpose  of classification.  A  highly used  strategy is
\emph{Dimension   Reduction}   (DR).    DR   aims  at   reducing   the
dimensionality of data  by mapping them onto another  space of a lower
dimension,  without discarding  any, or  as less  as possible,  of the
meaningful  information.    Recent  overviews  of  DR   can  be  found
in~\cite{DR:guided:tour,nl:dim:reduc,Guyon:2003:IVF:944919.944968}.
Two main approaches can be defined. 1) Unsupervised DR: The algorithms
are  applied  directly  on  the  data  without  exploiting  any  prior
information,  and project  the data  into a  lower  dimensional space,
according  to  some criterion  (data  variance  maximization for  PCA,
independence for ICA \ldots).   2) Supervised DR: Training samples are
available and are exploited to find a lower dimensional subspace where
the  class  separability is  improved.   Fisher Discriminant  Analysis
(FDA)  is  surely  one  of  the  most  famous  supervised  DR  method.
\textcolor{black}{However,  FDA maximizes  the ratio  of  the ``between
  classes'' scatter matrix, $S_b$,  and the ``within classes'' scatter
  matrix,  $S_w$.  The optimal  solution is  given by  the eigenvector
  corresponding  to the  first eigenvalues  of $S_w^{-1}S_b$.   In HD,
  $S_w^{-1}$   is  in   general  ill-conditioned   which   limits  the
  effectiveness  of the  method.}  Other  popular DR  methods  such as
Laplacian eigenmaps, Isomap  or Locally Linear Embedding~\cite[Chapter
4 and 5]{nl:dim:reduc} may be also limited by the dimensionality since
they are  based on  the Euclidean distance  between the  samples.  One
last  drawback  of   DR  methods  is  the  risk   of  losing  relevant
information.   In general,  DR methods  act globally,  which can  be a
problem for  classification purpose:  Different classes may  be mapped
onto the same subspace, even  if the global discrimination criteria is
maximized.

An alternative  strategy to DR  has been recently proposed,  i.e., the
subspace models~\cite{Parsons:2004:SCH:1007730.1007731}.  These models
assume that each class is located  in a specific subspace and consider
the original space  without DR for the processing.   For instance, the
Probabilistic Principal Component  Analysis (PPCA)~\cite{ppca} assumes
that  the classes  are  normally distributed  in  a lower  dimensional
subspace  and are  linearly  embedded in  the  original subspace  with
additive  white noise.   Such  models exploit  the \emph{empty  space}
property  of  HD   data  without  discarding  any   dimension  of  the
data~\cite{doi:10.1162/089976699300016728,5184847}.      A     general
subspace  model that  encompasses  several other  models  is the  High
Dimensional Discriminant Analysis (HDDA)  model, proposed by Bouveyron
\emph{et
  al.}~\cite{BOUVEYRON:2007:INRIA-00176283:1,BOUVEYRON:2007:HAL-00022183:4}.

Conversely, kernel based methods  do not reduce the dimensionality but
rather   work  with   the  full   HD   data~\cite{1151.30007}.   These
discriminative  methods are known  to be  more robust  to size  of the
dimensionality than  conventional generative methods.   However, local
kernel    methods    are    sensitive    to   the    size    of    the
dimensionality~\cite{NIPS2005_424}.   A kernel  method is  said  to be
local if the  decision function value for a new  sample depends on the
neighbors of  that sample in the  training set.  Since in  HD data the
neighborhood  of a  sample is  mostly  empty, such  local methods  are
negatively impacted by the dimension.  For instance, SVM with Gaussian
kernel
\begin{eqnarray}\label{eq:gaussian}
  k_g(\mathbf{x},\mathbf{z})=\exp\Bigg(-\frac{\|\mathbf{x}-\mathbf{z}\|^2}{2\sigma^2}\Bigg)
\end{eqnarray}
is such a local kernel method.

In this  paper, it is proposed  to use subspace models  to construct a
kernel  adapted to  high dimensional  data.  The  chosen  approach for
including  subspace models  in a  kernel function  is to  consider the
Mahalanobis distance, $d_{\boldsymbol{\Sigma}_c}$, between two samples
for    a     given    class,    $c$,     with    covariance    matrix,
$\boldsymbol{\Sigma}_c$:
$$d_{\boldsymbol{\Sigma}_c}(\mathbf{x},\mathbf{z}) = \sqrt{(\mathbf{x}-\mathbf{z})^t\boldsymbol{\Sigma}_c^{-1}(\mathbf{x}-\mathbf{z})}.$$
Previous          works           on          the          Mahalanobis
kernel~\cite{camps:mahalanobis,abe:training,abe:regression,4298136}        were
limited  by the  effect  of dimensionality  on  the matrix  inversion.
In~\cite{camps:mahalanobis}, the covariance matrix was computed on the
whole training set.  The associated implicit model is that the classes
share  the same  covariance matrix,  which  is not  true in  practice.
Diagonal    and   full    covariance   matrices    were   investigated
in~\cite{abe:training}   for  the   purpose   of  classification   and
in~\cite{abe:regression} for the purpose  of regression. However, in a
similar way, the  covariance matrix was computed for  all the training
samples. Computing the covariance  matrix for the Mahalanobis distance
with all the training samples is equivalent to project the data on all
the principal components, scale the variance to one, and then applying
the Euclidean distance.  By doing  so, classes could overlap more than
in the original input space  and the discrimination between them would
be decreased.

In this  work, the HDDA  model is used for  the definition of  a class
specific covariance  matrix adapted for  HD data. The  specific signal
and noise subspaces are estimated  for each considered class, ensuring
a  parsimonious characterization  of the  classes. Following  the HDDA
model it  is then possible  to derive  an explicit formulation  of the
inverse  of  the  covariance  matrix, without  any  regularization  or
dimension   reduction.   The   parsimonious   Mahalanobis  kernel   is
constructed   by  substituting   the  Euclidean   distance  with   the
Mahalanobis distance computed using the HDDA model.  It is proposed in
this work to  define several hyperparameters in the  kernel to control
the influence of the signal  and noise subspaces in the classification
process.   These hyperparameters  are  optimized  during the  training
process  by the  minimization  of the  so-called \emph{radius  margin}
bound of  the SVM classifier.  Compared  to the previous works  on the
Mahalanobis kernel for HD data, the  proposed method allows the use of
a more  complex model,  a separate covariance  matrix per  class, with
higher efficiency in terms of accuracy.  The remainder of the paper is
organized as follows.  The subspace  model and the proposed kernel are
discussed in  Section~\ref{sec:kernel}.  The problem of  selecting the
hyperparameters   for  classification   with  SVM   is  addressed   in
Section~\ref{sec:svm}. The Section~\ref{sec:pc} details the estimation
of the size of the signal subspace. Results on simulated and real high
dimensional data  are reported in  Section~\ref{sec:exp}.  Conclusions
and perspectives conclude the paper.

\begin{table}\small
  \centering
  \caption{Summary of HD spaces properties.}
  \label{table:hd:space}
  \begin{tabular}{|l|l|}
    \hline
    \multicolumn{2}{|c|}{\bf High Dimensional Spaces}\\
    \hline
    \multicolumn{1}{|c|}{\em Curse} & \multicolumn{1}{c|}{\em Blessing}\\
    \hline
    Poor statistical estimation & Emptiness\\
    Concentration of measure  & Class separability \\
    \hline
  \end{tabular}
\end{table}

\section{Regularized Mahalanobis Kernel}\label{sec:kernel}
\subsection{Review of HDDA model}
The most general  HDDA sub-model is used  in this work\footnote{Refers
  to                                               $[a_{ij}b_iQ_id_i]$
  in~\cite{BOUVEYRON:2007:INRIA-00176283:1,BOUVEYRON:2007:HAL-00022183:4}.},
i.e., each class has his own  specific subspace.  Here, we will review
the HDDA model but restricted to  the problem of the covariance matrix
inversion. However HDDA was  originally proposed for classification or
clustering with Gaussian mixture model.  Interested readers can find a
detailed               presentation              of               HDDA
in~\cite{BOUVEYRON:2007:INRIA-00176283:1,BOUVEYRON:2007:HAL-00022183:4}.

\begin{figure}
  \centering \resizebox{0.49\textwidth}{!}{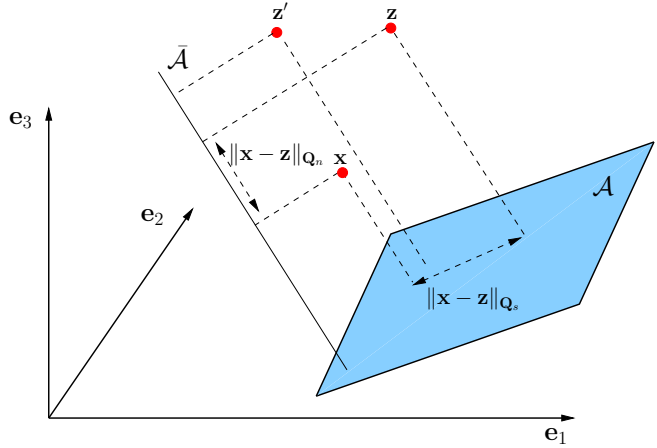}
  \caption{Cluster-based model. The  distance between $\mathbf{x}$ and
    $\mathbf{z}$ is  computed both in  the signal subspace and  in the
    noise     subspace.      Note     that     in     this     example
    $\dim(\bar{\mathcal{A}}_c)<\dim(\mathcal{A}_c)$, but for real data
    it is usually the opposite.  $\|\cdot\|_{\mathbf{Q}_s}$ is the dot
    product in  $\mathcal{A}_c$ and $\|\cdot\|_{\mathbf{Q}_n}$  is the
    dot product in $\bar{\mathcal{A}}_c$.}
  \label{fig:ppca}
\end{figure}

In subspace  models, it is assumed  that the data from  each class are
clustered in the vector space.  This  cluster does not need to have an
elliptic  shape but it  is generally  assumed that  the data  follow a
Gaussian distribution.  The covariance matrix  of the class $c$ can be
written through its eigenvalue decomposition:
$$\boldsymbol{\Sigma}_c = \mathbf{Q}_c\boldsymbol{\Lambda}_c\mathbf{Q}_c^t$$
where $\boldsymbol{\Lambda}_c$  is the diagonal  matrix of eigenvalues
$\lambda_{ci}$,  $i\in\{1,\ldots,d\}$, of  $\boldsymbol{\Sigma}_c$ and
$\mathbf{Q}_c$   is  the  matrix   that  contains   the  corresponding
eigenvectors~$\mathbf{q}_{ci}$.   The  HDDA  model assumes  the  $p_c$
first eigenvalues are different  and the remaining $d-p_c$ eigenvalues
are identical. The  model is similar to PPCA, but  more general in the
sense that  additional sub-models can  be defined. In  particular, the
intrinsic dimension  $p_c$ are not  constrained in HDDA  whereas there
are  assumed to  be equal  for each  class in  PPCA.\\ Under  the HDDA
framework, the covariance matrix has the following expression:
$$\boldsymbol{\Sigma}_c=\sum_{i=1}^{p_c}\lambda_{ci}\mathbf{q}_{ci}\mathbf{q}_{ci}^t + b_c\sum_{i={p_c}+1}^d\mathbf{q}_{ci}\mathbf{q}_{ci}^t$$
where the last $d-p_c$ eigenvalue are equal to $b_c$.  The inverse can
be computed explicitly by
$$\boldsymbol{\Sigma}_c^{-1}=\underbrace{\sum_{i=1}^{p_c}\frac{1}{\lambda_{ci}}\mathbf{q}_{ci}\mathbf{q}_{ci}^t}_{\mathcal{A}_c} + \underbrace{\frac{1}{b_c}\sum_{i={p_c}+1}^d\mathbf{q}_{ci}\mathbf{q}_{ci}^t}_{\bar{\mathcal{A}}_c}.$$
This statistical model can be understood equivalently by a geometrical
assumption: For each class, the data belong to a cluster that lives in
a lower dimensional space $\mathcal{A}_c$, namely the signal subspace.
The     original    input     space    can     be     decomposed    as
$\mathbb{R}^d=\mathcal{A}_c\bigoplus\bar{\mathcal{A}}_c$            (by
construction $\bar{\mathcal{A}}$ is  the noise subspace which contains
only  white noise).   Figure~\ref{fig:ppca} gives  an  illustration of
that in~$\mathbb{R}^3$.

Using    $\mathbf{I}   =\sum_{i=1}^d\mathbf{q}_{ci}\mathbf{q}_{ci}^t$,
$\mathbf{I}$  being the identity  matrix, the  inverse can  be finally
written as
\begin{eqnarray}\label{eq:inv}\boldsymbol{\Sigma}_c^{-1}=\sum_{i=1}^{p_c}\Big(\frac{1}{\lambda_{ci}}-\frac{1}{b_c}\Big)\mathbf{q}_{ci}\mathbf{q}_{ci}^t + \frac{1}{b_c}\mathbf{I}.\end{eqnarray}
Standard   likelihood   maximization   shows   that   the   parameters
$(\lambda_{ci},\mathbf{q}_{ci})_{i=1,\dotsc,p_c}$  and  $b_c$  can  be
computed          from          the         sample          covariance
matrix~\cite{BOUVEYRON:2007:INRIA-00176283:1}:
$$\hat{\boldsymbol{\Sigma}}_c=\frac{1}{n_c}\sum_{i=1}^{n_c}\big(\mathbf{x}_i-\bar{\mathbf{x}}_c\big)\big(\mathbf{x}_i-\bar{\mathbf{x}}_c\big)^t$$
where $\bar{\mathbf{x}}_c$ is the sample  mean for class $c$ and $n_c$
the number  of samples of  the class.  $\lambda_{ci}$ is  estimated by
the    $i$   first   eigenvalue    of   $\hat{\boldsymbol{\Sigma}}_c$,
$\mathbf{q}_{ci}$  by the corresponding  eigenvector and  $\hat{b}_c =
\big(\text{trace}(\hat{\boldsymbol{\Sigma}}_c)-\sum_{i=1}^{\hat{p}_c}\hat{\lambda}_{ci}\big)/(d-\hat{p}_c)$
(the estimation  of the dimension  $p_c$ of the subspace  is discussed
later).   The   last  $d-p_c$  eigenvalues   and  their  corresponding
eigenvectors  are  not  needed  for  the computation  of  the  inverse
in~(\ref{eq:inv}).

The major advantage of such a model is that it reduces drastically the
number  of parameters to  estimate for  computing the  inverse matrix.
Indeed, with the full  covariance matrix, $d(d+3)/2$ parameters are to
be estimated.   With the HDDA model, only  $d(p_c+1) +1 -p_c(p_c-1)/2$
parameters are to be estimated. For instance, if $d=100$ and $p_c=10$,
5150 parameters are  needed for the full covariance  and only 1056 for
the  HDDA model.   Furthermore, the  stability is  improved  since the
smallest eigenvalues of the  covariance matrix and their corresponding
eigenvectors, which are difficult  to compute accurately, are not used
in~(\ref{eq:inv}).\\
Finally,  using the HDDA  model, the  square Mahalanobis  distance for
class $c$ is approximated by
\begin{eqnarray}\label{eq:dist:maha}
  d^2_{\hat{\boldsymbol{\Sigma}}_c}(\mathbf{x},\mathbf{z}) = \sum_{i=1}^{\hat{p}_c}\Big(\frac{1}{\hat{\lambda}_{ci}}-\frac{1}{\hat{b}_c}\Big)\|\hat{\mathbf{q}}_{ci}^t(\mathbf{x}-\mathbf{z})\|^2 + \frac{\|\mathbf{x}-\mathbf{z}\|^2}{\hat{b}_c}.
\end{eqnarray}

This  approach relies  on  the analysis  of  the empirical  covariance
matrix,  as  with  PCA.   But  instead  of  keeping  only  significant
eigenvalues,  (\ref{eq:dist:maha}) considers  all the  original space,
without  discarding  any dimension.   This  has  two main  theoretical
advantages over the conventional PCA:

\begin{enumerate}
\item Two samples may be close in the signal subspace but far apart in
  the original space, which is a problem for classification tasks.  It
  can be handled  by considering the noise subspace  together with the
  signal subspace.  Consider  for instance $\mathbf{z}$, $\mathbf{z}'$
  and   $\mathbf{x}$  in  Figure~\ref{fig:ppca}.    In  $\mathcal{A}$,
  $\mathbf{z}'$ seems closer  to $\mathbf{x}$ than $\mathbf{z}$, while
  it is  not as it  can be seen  by adding $\bar{\mathcal{A}}$  in the
  distance computation.
\item An  accurate estimation of the signal  subspace size $\hat{p}_c$
  is  necessary   if  PCA  is   applied:  The  worst   scenario  being
  $\hat{p}_c<<p_c$, \emph{i.e.},  relevant eigenvectors are discarded.
  By considering both  the signal and the noise  subspaces, the method
  becomes  less sensitive  to  $\hat{p}_c$.  Even  in  the worst  case
  scenario, the eigenvectors are still considered.
\end{enumerate}

\subsection{Mahalanobis Kernel}
\begin{figure*}
  \centering
  \begin{tabular}{ccc}
    \includegraphics[width=0.3\textwidth]{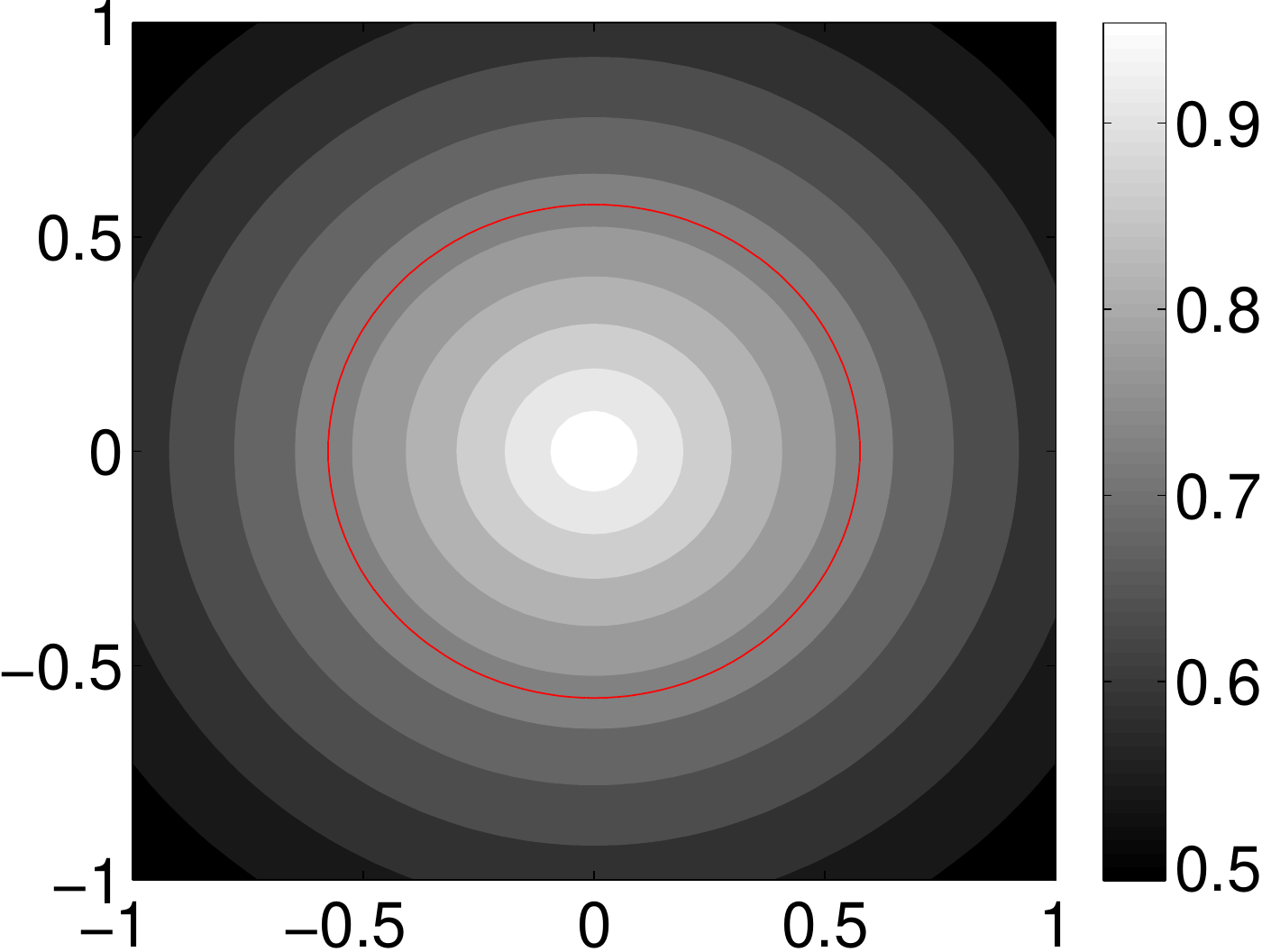} &  \includegraphics[width=0.3\textwidth]{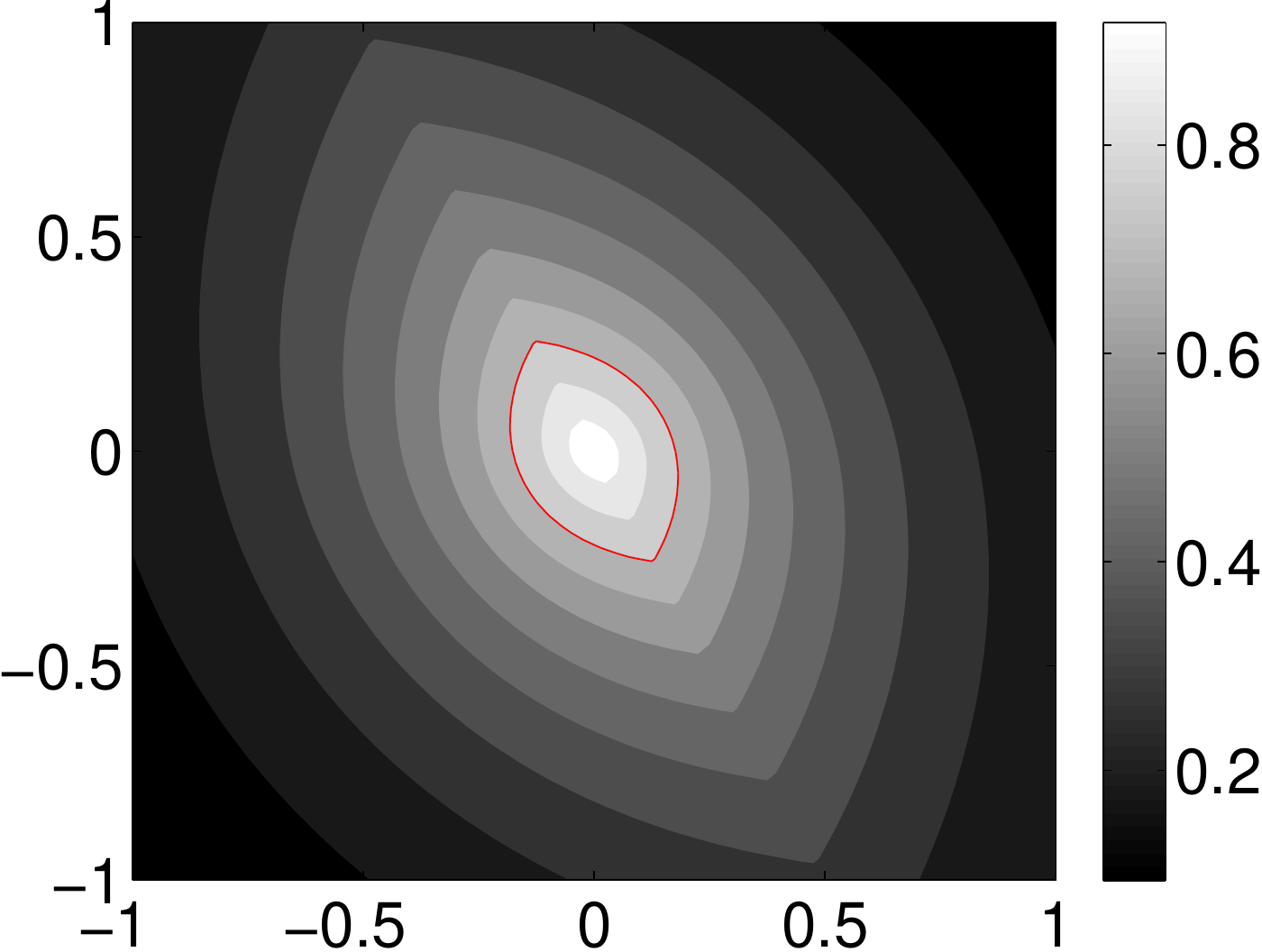} &  \includegraphics[width=0.3\textwidth]{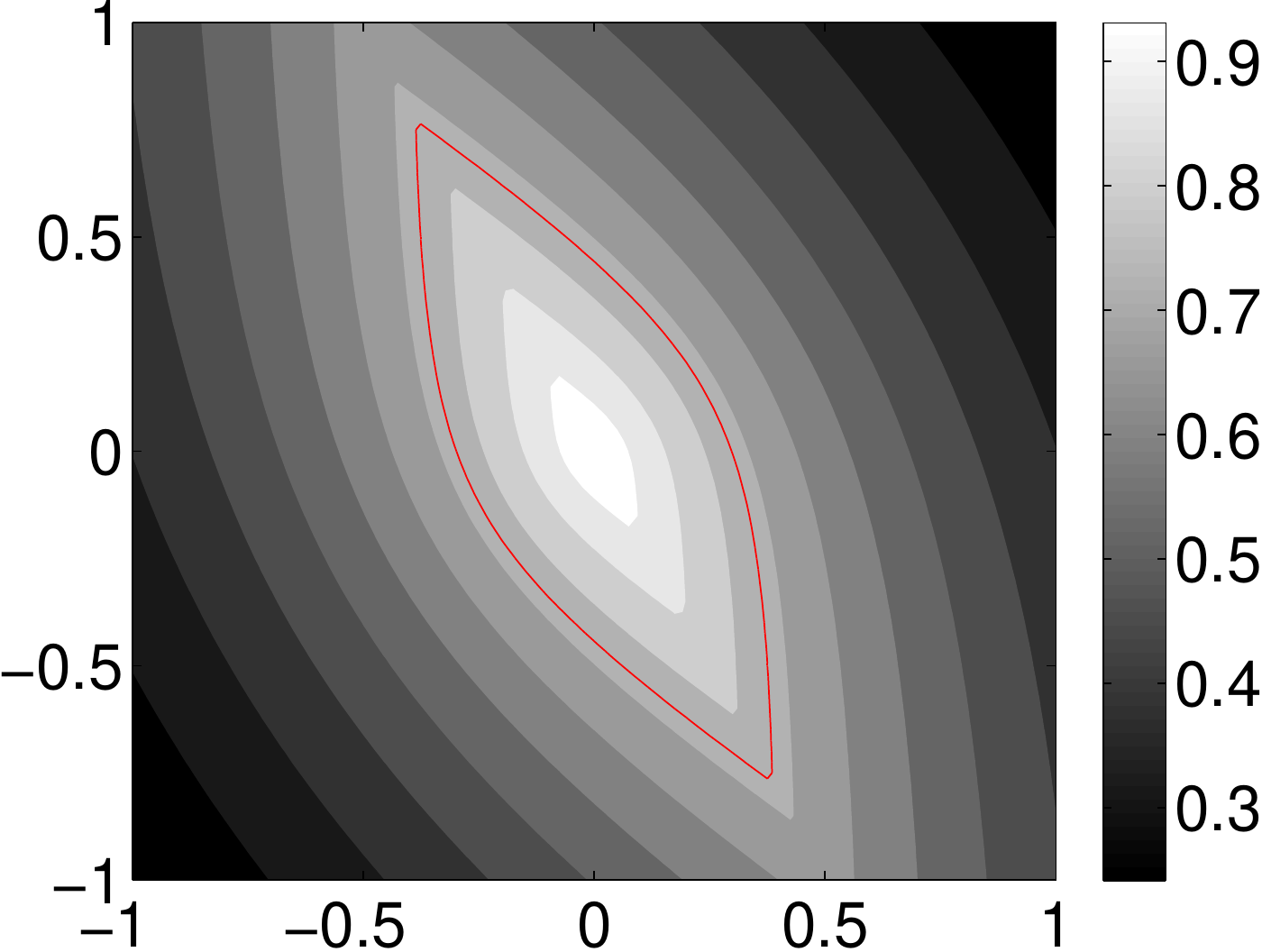}\\
    (a) & (b) & (c)
  \end{tabular}
  \caption{Values  of the  kernel  function $k(\mathbf{0},\mathbf{x})$
    with $\mathbf{0}=[0,0]$ and $\mathbf{x}\in[-1,1]\times[-1,1]$. The
    red line represents  the contour line for the  value 0.75. (a) is
    the  Gaussian kernel,  (b)  is the  kernel (\ref{eq:kernel})  with
    $\sigma_1^2   =   \sigma_2^2   =   0.5$,   (c)   is   the   kernel
    (\ref{eq:kernel})  with   $\sigma_1^2  =1.5$  and   $\sigma_2^2  =
    0.5$. The covariance matrix used was $[0.6~-0.2;-0.2~0.6]$ and the
    signal  subspace  was  of   dimension  1,  spanned  by  the  first
    eigenvector of the covariance matrix.}
  \label{fig:kernel}
\end{figure*}
\textcolor{black}{The regularized  Mahalanobis kernel for  class $c$ is
  constructed  by substituting  (\ref{eq:dist:maha}) to  the Euclidean
  distance  in the  Gaussian kernel~(\ref{eq:gaussian})  and switching
  eigenvalues   $(\hat{\lambda}_{ci},\hat{b}_c)$   to  hyperparameters
  $(\sigma_{ci}^2,\sigma_{c\hat{p}_c+1}^2)$ that  are optimized during
  the training step:}
\begin{eqnarray}\label{eq:kernel}
  \begin{array}{l}
    k_{m}(\mathbf{x},\mathbf{z}|c) =\\ \displaystyle{\exp\Bigg(-\frac{1}{2}\bigg( \sum_{i=1}^{\hat{p}_c}\frac{\|\hat{\mathbf{q}}_{ci}^t(\mathbf{x}-\mathbf{z})\|^2}{\sigma_{ci}^2} + \frac{\|\mathbf{x}-\mathbf{z}\|^2}{\sigma_{c\hat{p}_c+1}^2}\bigg)\Bigg)}
  \end{array}
\end{eqnarray}
where    $\sigma_{ci},\    i\in\{1,\dotsc,\hat{p}_c+1\}$    are    the
hyperparameters of the  kernel. As described in Section~\ref{sec:svm},
these   parameters  are   tuned   during  the   training  step.    The
hyperparameters have been introduced  for the following reason.  It is
known that the principal directions are not optimal for classification
since  they do  not maximize  any discrimination  criterion.  However,
they still span  a subspace where there are variations  in the data of
the  considered  class.  The  hyperparameters  $\sigma_{ci}$ allow  to
control which directions are more relevant (or discriminative) for the
classification  process:  The feature  space  is  modified during  the
training process to ensure a better discrimination between samples.

It is interesting to note  that the regularized Mahalanobis kernel can
be expressed as the product of Gaussian kernels:
\begin{eqnarray}\label{eq:mixture:prod}
  k_m(\mathbf{x},\mathbf{z}|c) = k_g(\mathbf{x},\mathbf{z})\times\prod_{i=1}^{\hat{p}_c} k_g(\hat{\mathbf{q}}_{ci}^t\mathbf{x},\hat{\mathbf{q}}_{ci}^t\mathbf{z}).
\end{eqnarray}
The  feature space  induced by  the kernel  and the  influence  of the
hyperparameters is analyzed in the next section.

\subsection{Geometry of the induced feature space}
Working  with a  kernel function  is equivalent  to work  with samples
mapped onto  a feature space  $\mathcal{H}$, where the dot  product is
equivalent    to    the     kernel    evaluation    in    the    input
space~\cite{1151.30007,Filippone2008176}:
$$k(\mathbf{x},\mathbf{z})=\langle \phi(\mathbf{x}), \phi(\mathbf{z})\rangle_{\mathcal{H}},$$
$\phi$  being  the  feature  map.   Under some  weak  conditions,  the
projected  samples   in  the  feature  space  live   on  a  Riemannian
manifold~\cite{kernel:space,willams:li:feng:wu}.  The metric tensor~is
\begin{eqnarray}
  g_{ij}(\mathbf{x})=\frac{\partial^2k(\mathbf{x},\mathbf{z})}{\partial x_i\partial z_j}\Bigg|_{\mathbf{z}=\mathbf{x}}
\end{eqnarray}
which        is,         for        the        Gaussian        kernel,
$g_{ij}(\mathbf{x})=\sigma^{-2}\delta_{ij}$  with  $\delta_{ij}=1$  if
$i=j$  and $0$  otherwise.  This  metric stretches  or  compresses the
Euclidean distance  between $\mathbf{x}$ and $\mathbf{z}$  by a factor
$\sigma^{-2}$. Each variable is assumed equally relevant for the given
task, \emph{e.g.}, classification or regression.

For the kernel~(\ref{eq:kernel}) the metric tensor is:
\begin{eqnarray}
  g_{ij}(\mathbf{x}|c) = \sum_{l=1}^{\hat{p}_c}\frac{q_{cli}q_{clj}}{\sigma_{cl}^2} + \frac{\delta_{ij}}{\sigma_{c\hat{p}_c+1}^2}
\end{eqnarray}
with   $q_{cli}$   the    $i^{th}$   element   of   $\mathbf{q}_{cl}$.
\textcolor{black}{The  distance between  two samples  is  stretched (if
  $\sigma_{cl}^2\geq  1$)  or  compressed (if  $\sigma_{cl}^2\leq  1$)
  along the $\hat{p}_c$ first principal components of class $c$ (first
  term  of the  right part  of the  equation) and  along  the original
  components  (last  term of  the  equation).   In  other words,  each
  principal component  is weighted according to its  relevance for the
  processing.}

The analysis of the metric tensor  exhibits the nature of the proposed
kernel for a given class: It is  a mixture of a Gaussian kernel on the
original  variables and  a Gaussian  kernel on  the $\hat{p}_c$  first
principal  components of  the  considered  class. The  hyperparameters
$\sigma_{cl}$ are tuned during the  training process.  This allows the
optimization of the weight of each kernel.  If $\sigma_{cl}^2=+\infty,
\forall l\in\{1,\dotsc,\hat{p}_c\}$, (\ref{eq:kernel})  reduces to the
conventional    Gaussian     kernel.     On    the     contrary,    if
$\sigma_{c\hat{p}_c+1}^2=+\infty$,  (\ref{eq:kernel})  reduces to  the
Gaussian  kernel  on  the   $\hat{p}_c$  first  principal  components.
Figure~\ref{fig:kernel} shows  the kernel values for  different values
of the hyperparameters. Note that opposite to the Gaussian kernel, the
kernel in~(\ref{eq:kernel}) is not isotropic.

The  following  section  reviews  the  basics of  SVM  classifier  and
presents how the hyperparameters are computed.
\section{L2-SVM and Radius margin bound optimization}\label{sec:svm}
\textcolor{black}{Support   vector  machines   (SVM)   is  a   standard
  classification kernel methods~\cite{Burges98atutorial}. It has shown
  to performs very  well on several data sets  from moderate dimension
  to                          high                         dimensional
  data~\cite{Kim20102871,Fauvel:2012:SKA:2030819.2031227}.    In   the
  following  section, the  main results  are presented  but interested
  readers                           could                          see
  references~\cite{Burges98atutorial,intro:svm:2000,Vapnik:2}       for
  further mathematical details about the SVM framework.}

\subsection{L2-Support Vector Machines}
The L2-SVM  is considered  in this work  rather than  the conventional
L1-SVM~\cite{intro:svm:2000}: With L2-SVM  it is possible to
tune  the hyperparameters  automatically by  optimizing the  so called
radius-margin   bound~\cite{chapelle:ams}.   The  L2-SVM   solves  the
conventional L1-SVM optimization problem with a quadratic penalization
of    errors~\cite{intro:svm:2000}.      Given    a    training    set
$\mathcal{S}=\{(\mathbf{x}_1,y_1),\dotsc,(\mathbf{x}_n,y_n)\},\
(\mathbf{x}_i,y_i)\in   \mathbb{R}^d\times\{-1;1\}$,   the  parameters
$(\alpha_i)_{i=1}^n$ and $b$ of the decision function $f$,
$$f(\mathbf{z})= \sum_{i=1}^n\alpha_ik(\mathbf{x}_i,\mathbf{z}) + b,$$ 
are found by solving the convex optimization problem:
\begin{eqnarray}\label{eq:svmtrain}
  \begin{array}{rcl}
    \underset{\alpha}{\max}\  g(\alpha) & = &\displaystyle{ \sum_{i=1}^n\alpha_i - \frac{1}{2}\sum_{i,j=1}^n\alpha_i\alpha_jy_iy_j \tilde{k}(\mathbf{x}_i,\mathbf{x}_j) }\\
    \text{subject to}          &   & 0\leq \alpha_i \text{ and } \displaystyle{\sum_{i=1}^n\alpha_iy_i =0}
  \end{array}
\end{eqnarray}
where                             $\tilde{k}(\mathbf{x}_i,\mathbf{x}_j)
=k(\mathbf{x}_i,\mathbf{x}_j)  +C^{-1}\delta_{ij}$  with  $k$  the
kernel  function and  $C$ a  positive hyperparameter  that is  used to
penalize the training errors.

An estimate of the generalization errors is given by an upper bound on
the number of errors of the leave-one-out procedure, the radius-margin
bound $\mathcal{T}$~\cite{Vapnik:2}:
\begin{eqnarray}
  \mathcal{T}(\mathbf{p}):= \mathcal{R}^2\mathcal{M}^2.
\end{eqnarray}
$\mathcal{R}^2$  is  the  radius  of  the  smallest  hypersphere  that
contains  all $\phi(\mathbf{x}_i)$, $\mathcal{M}^2$  is the  margin of
the  classifier, it  is given  by  the optimal  objective function  of
(\ref{eq:svmtrain}), and $\mathbf{p}$  are the hyperparameters. In our
setting      $\mathbf{p}=[\sigma^2_1,\dotsc,\sigma^2_{\hat{p}_c+1},C]$.
$\mathcal{R}^2$ is  obtained by the optimal objective  function of the
following constraint optimization problem~\cite{Vapnik:2}:
\begin{eqnarray}\label{eq:radius}
  \begin{array}{rcl}
    \underset{\beta}{\max}\   g'(\beta)&=&\displaystyle{\sum_{i=1}^n\beta_i \tilde{k}(\mathbf{x}_i,\mathbf{x}_i)-\sum_{i,j=1}^n\beta_i\beta_j \tilde{k}(\mathbf{x}_i,\mathbf{x}_j)}\\
    \text{subject to}& & 0\leq \beta_i \text{ and } \displaystyle{\sum_{i=1}^n\beta_i=1}.
\end{array}
\end{eqnarray}
Since both $\mathcal{R}^2$ and $\mathcal{M}^2$ depend on $\mathbf{p}$,
it is  possible to optimize $\mathcal{T}$ to  set the hyperparameters.
Chapelle \emph{et  al.}~\cite{chapelle:ams}, followed later  by S.  S.
Keerthi~\cite{1031955}, have  proposed an algorithm  based on gradient
optimization method.  It is discussed in the following section.

With the proposed  kernel, if two classes $i$ and  $j$ are considered,
the classifier  for ``$i$ vs $j$''  is not the same  as the classifier
for  ``$j$ vs  $i$''  since the  kernel function  is  specific to  the
classes $i$ and  $j$, respectively. Indeed, for  a multiclass problem,
the ``one vs  one'' approach must not  be used and the  ``one vs all''
approach should be preferred~\cite{991427}.

\subsection{Radius-margin bound Optimization}
Computing the  gradient of  $\mathcal{T}$ requires the  computation of
the gradient of the following expressions\footnote{For simplicity, the
  parameter   $c$   of  the   kernel   function   is  omitted,   i.e.,
  $\tilde{k}(\mathbf{x},\mathbf{z}|c)$       is       written       as
  $\tilde{k}(\mathbf{x},\mathbf{z})$.}
\begin{eqnarray}\label{eq:w2}
  \mathcal{M}^2& =& 2 \sum_{i=1}^n \tilde{\alpha}_i - \sum_{i,j=1}^n\tilde{\alpha}_i\tilde{\alpha}_jy_iy_j \tilde{k}(\mathbf{x}_i,\mathbf{x}_j)
\end{eqnarray}
and of
\begin{eqnarray}\label{eq:r2}
\mathcal{R}^2 &= &\sum_{i=1}^n\tilde{\beta}_i \tilde{k}(\mathbf{x}_i,\mathbf{x}_i)-\sum_{i,j=1}^n\tilde{\beta}_i\tilde{\beta}_j \tilde{k}(\mathbf{x}_i,\mathbf{x}_j)
\end{eqnarray}
where $(\tilde{\alpha}_i)_{i=1}^n$ and $(\tilde{\beta}_i)_{i=1}^n$ are
the optimal  parameters of (\ref{eq:svmtrain})  and (\ref{eq:radius}).
The  gradient of  (\ref{eq:w2}) depends  on  $\tilde{\alpha}_i$, which
depends  on $\mathbf{p}$  (similar comments  hold  for (\ref{eq:r2})).
Chapelle  \emph{et al.}   have proven  that since  $\mathcal{M}^2$ and
$\mathcal{R}^2$  are   computed  via  an   optimization  problem,  the
gradients  of $\tilde{\alpha}_i$  and $\tilde{\beta}_i$  do  not enter
into      account      in      the      computation      of      their
gradients~\cite{chapelle:ams}.      Hence,     the     gradient     of
(\ref{eq:radius}) can be written as:
\begin{eqnarray}
  \nabla\mathcal{T}=\Bigg[\frac{\partial \mathcal{T}}{\partial C},\frac{\partial \mathcal{T}}{\partial \sigma_1^2},\ldots,\frac{\partial \mathcal{T}}{\partial \sigma_{\hat{p}_c+1}^2}\Bigg]^t
\end{eqnarray}
with 
\begin{eqnarray}
  \frac{\partial \mathcal{T}}{\partial C}=\frac{\partial \mathcal{R}^2}{\partial C}\mathcal{M}^2 + \mathcal{R}^2\frac{\partial \mathcal{M}^2}{\partial C}
\end{eqnarray}
and 
\begin{eqnarray}
  \frac{\partial \mathcal{T}}{\partial \sigma^2_\ell}=\frac{\partial \mathcal{R}^2}{\partial \sigma_\ell^2}\mathcal{M}^2 + \mathcal{R}^2\frac{\partial \mathcal{M}^2}{\partial \sigma^2_\ell}
\end{eqnarray}
for $\ell\in\{1,\dots,\hat{p}_c+1\}$. The derivatives of $\mathcal{R}^2$ are
\begin{eqnarray}
  \frac{\partial \mathcal{R}^2}{\partial C} &= &\frac{1}{C^2}\sum_{i=1}^n\tilde{\beta}_i(\tilde{\beta}_i -1)\\
  \frac{\partial \mathcal{R}^2}{\partial \sigma_\ell^2} &=& -\sum_{i,j=1}^n\tilde{\beta}_i\tilde{\beta}_j\frac{\partial \tilde{k}(\mathbf{x}_i,\mathbf{x}_j)}{\partial \sigma_\ell^2}
\end{eqnarray}
with 
\begin{eqnarray}
  \begin{array}{l}
  \displaystyle{\frac{\partial \tilde{k}(\mathbf{x}_i,\mathbf{x}_j)}{\partial \sigma_\ell^2} =} \\
  \left\{
    \begin{array}{ll}
      \displaystyle{\frac{\|\hat{\mathbf{q}}_{c\ell}^t(\mathbf{x}_i-\mathbf{x}_j) \|^2}{\sigma^4_\ell}}\tilde{k}(\mathbf{x}_i,\mathbf{x}_j) & \text{if } \ell \in\{1,\ldots,\hat{p}_c\}\\
      \displaystyle{\frac{\|\mathbf{x}_i-\mathbf{x}_j\|^2}{\sigma^4_\ell}}\tilde{k}(\mathbf{x}_i,\mathbf{x}_j) & \text{if } \ell = \hat{p}_c+1.
    \end{array}\right.
  \end{array}
\end{eqnarray}
The derivatives of $\mathcal{M}^2$ are
\begin{eqnarray}
  \frac{\partial \mathcal{M}^2}{\partial C} &= &\frac{1}{C^2}\sum_{i=1}^n\tilde{\alpha}_i\\  
  \frac{\partial \mathcal{M}^2}{\partial \sigma_\ell^2} &=& -\sum_{i,j=1}^n\tilde{\alpha}_i\tilde{\alpha}_jy_iy_j\frac{\partial \tilde{k}(\mathbf{x}_i,\mathbf{x}_j)}{\partial \sigma_\ell^2}.
\end{eqnarray}
\textcolor{black}{Once   the  derivatives   have  been   computed,  the
  optimization  of  $\mathcal{T}$   is  done  through  a  conventional
  gradient descent, following the framework in~\cite{chapelle:ams}. At
  each iteration $t$, the set of hyperparameters are updated as with a
  step proportional to the negative of the gradient of $\mathcal{T}$:
  $$\mathbf{p}^{t+1} = \mathbf{p}^{t} - \gamma\nabla\mathcal{T}$$
  where $\gamma\geq  0$ is a  step size parameter.  For implementation
  details, see~\cite{ams:online}}.

\section{Estimation of $\hat{p}_c$}\label{sec:pc}
\begin{figure}
  \centering
  \includegraphics[width=0.3\textwidth]{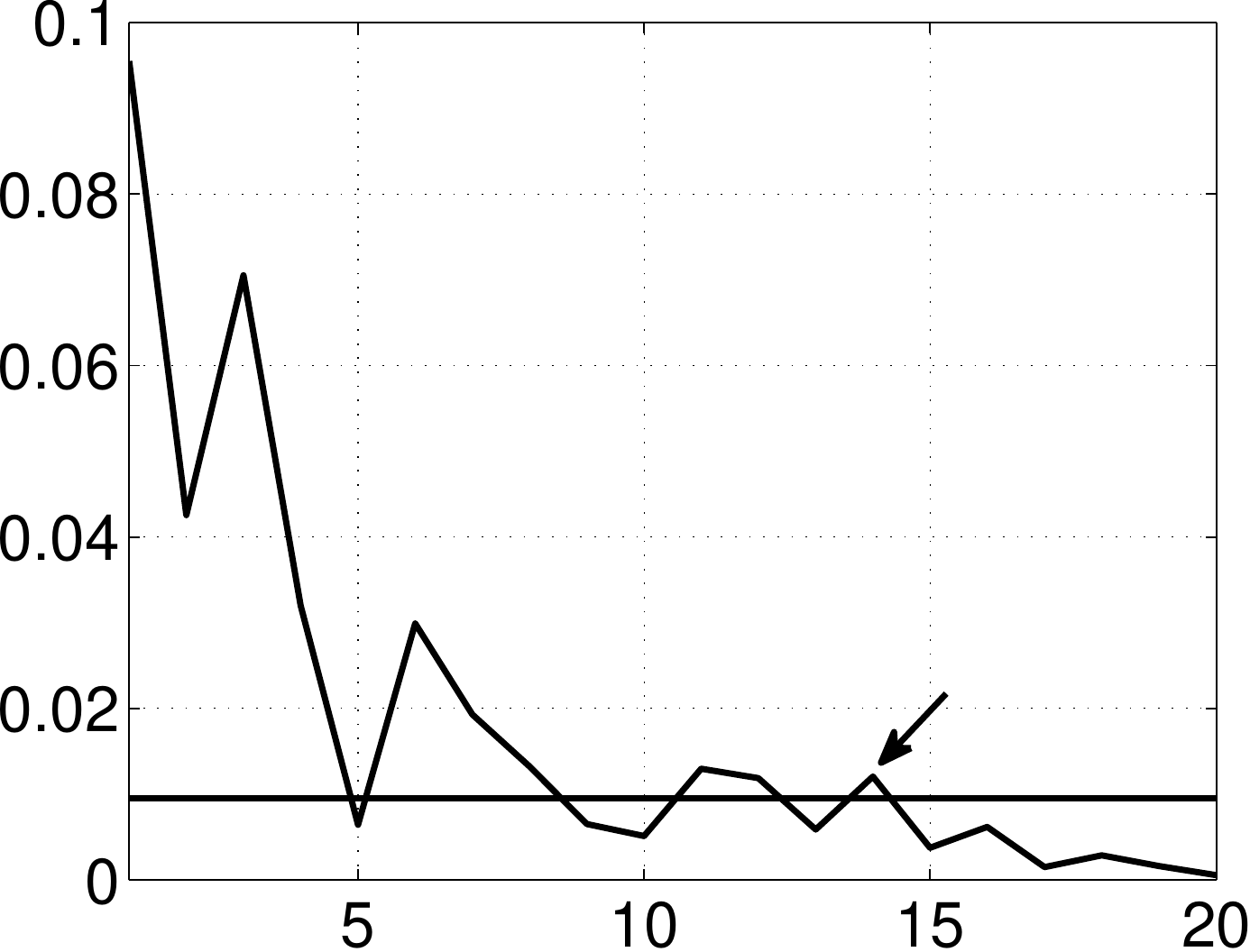}
  \caption{\textcolor{black}{Scree test of  Cattell.  The threshold $s$
      is  set  to  10\%.   The   estimated  $\hat{p}_c$  is  14.   See
      Section~\ref{sec:simu} for  a description of the  data set.  The
      horizontal  axis represents  the index  $i$.  The  vertical axis
      represents  the  numerical  difference between  two  consecutive
      eigenvalues,  $\Delta_i$.  The  curve represents  the difference
      between  two consecutive  eigenvalues  and  the horizontal  line
      represents 10 \% of the highest difference.  The arrow shows the
      estimated $\hat{p}_c$.}}
  \label{fig:cattell}
\end{figure}

The size  of the signal  subspace was estimated  by the scree  test of
Cattell~\cite{citeulike:3574985}   using  the   same   methodology  as
in~\cite{BOUVEYRON:2007:INRIA-00176283:1}.    The  test   consists  in
comparing   the  difference,   $\Delta_i$,  between   two  consecutive
eigenvalues          $\lambda_i$          and         $\lambda_{i+1}$,
$\Delta_i=\lambda_i-\lambda_{i+1}$.   When the  differences $\Delta_i$
are   below  a  user-defined   threshold  $s$   for  all   $i$,  i.e.,
$\Delta_j<s,\forall  j \in\{i,\ldots,d-1\}$, $\hat{p}_c$  is estimated
as $\hat{p}_c=i$.   In general  the threshold is  a percentage  of the
highest  difference.  Figure~\ref{fig:cattell} shows  an example  on a
simulated data  set (see  Section~\ref{sec:simu} for a  description of
the data).  The correct value in  that case is $p=10$ but its estimate
is $\hat{p}=14$.

\section{Experimental results}\label{sec:exp}
Classification results  are presented in this  section.  Regarding the
multiclass  strategy, the  results  must be  considered as  individual
binary  classification problems:  No fusion  rules were  applied.  For
instance,   in  Table~\ref{tab:res},   the  results   for  the   class
``Asphalt'' should be read as ``Asphalt  vs all''.  The reason for the
use of this approach is the better interpretation of the results which
is  obtained because  the results  are  not biased  by the  multiclass
fusion strategy.

\begin{figure*}
  \centering
  \begin{tabular}{ccc}
    \includegraphics[width=0.3\textwidth]{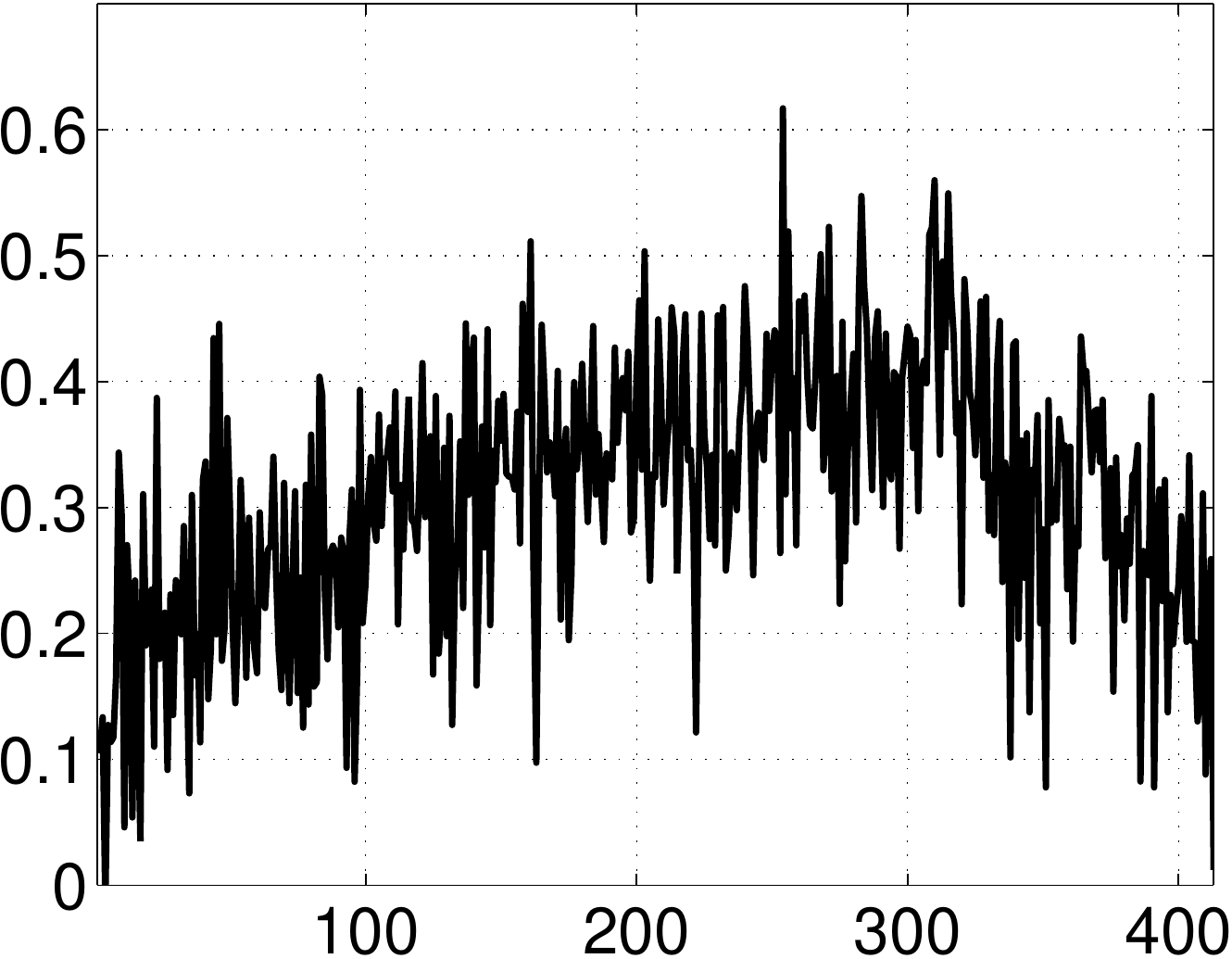} & \includegraphics[width=0.3\textwidth]{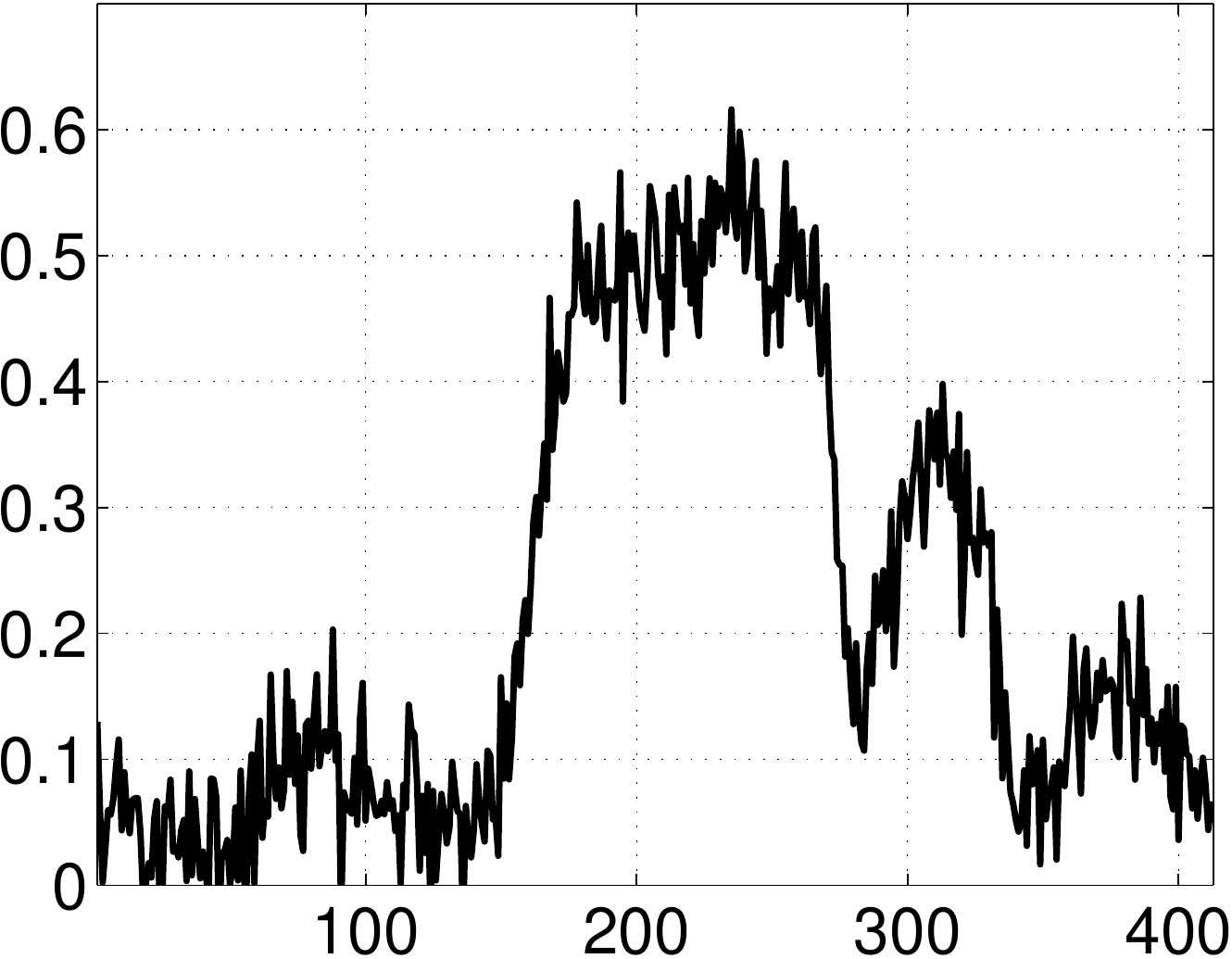} & \includegraphics[width=0.3\textwidth]{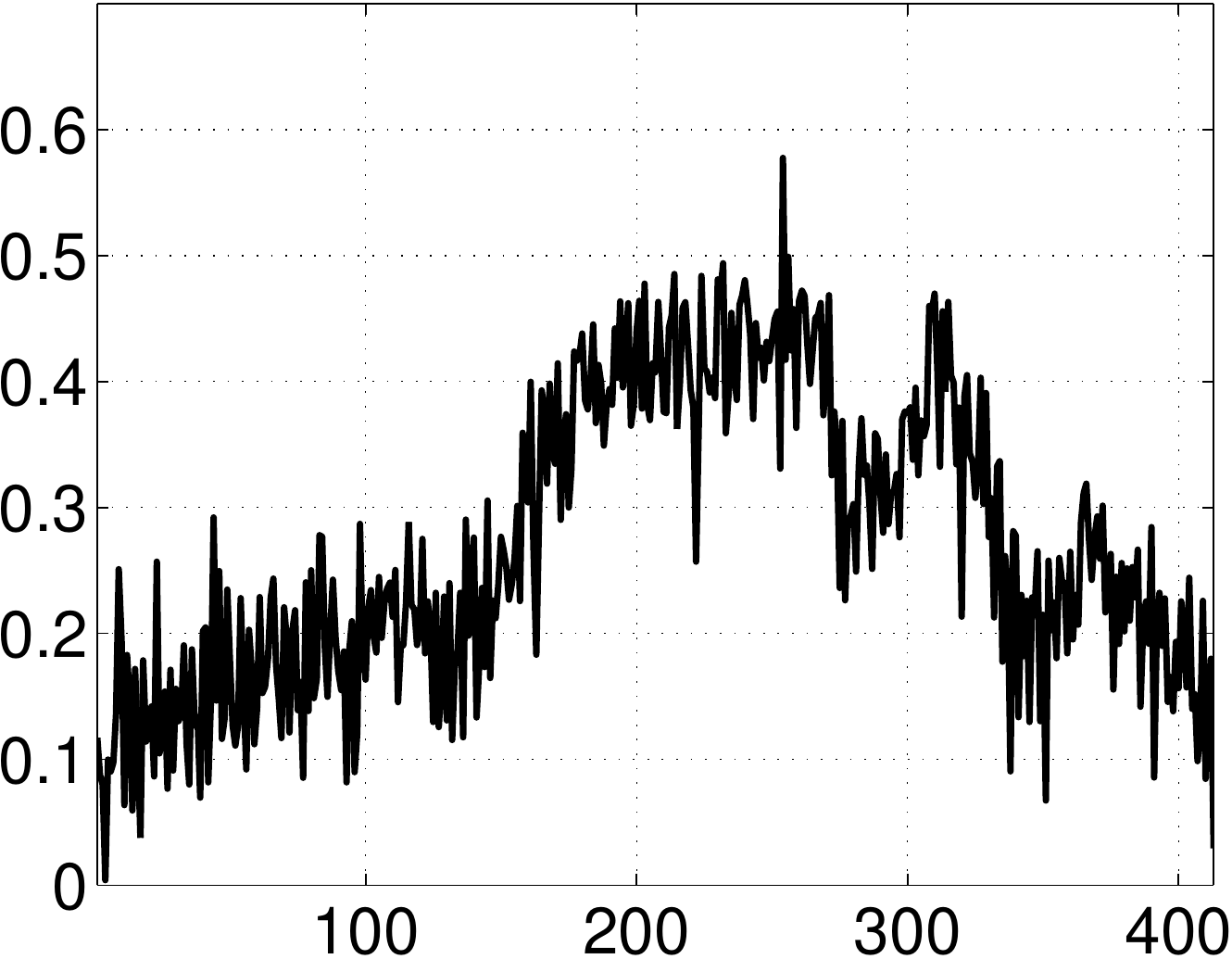}\\
    (a) & (b) & (c)
  \end{tabular}
  \caption{Simulated  spectra: (a) $\mathbf{s}_1$,  (b) $\mathbf{s}_2$
    and (c) $\mathbf{x}$.  The horizontal axis is the variable and the
    vertical axis is the simulated reflectance.  The parameters of the
    simulation  are: $\alpha=[0.6,0.4]$,  $d=413$,  $p=10$, $N_c=2$  and
    SNR=1.}
  \label{fig:simulated:spectra}
\end{figure*}
\subsection{Classification   of    simulated   data   following   HDDA
  model}\label{sec:simu}
\textcolor{black}{In  this section,  the  proposed  kernel, namely  the
  HDDA-Mahalanobis  Kernel   (HDDA-MK),  is  used  with   the  SVM  in
  classification and evaluated on simulated data.  The performances in
  terms of classification accuracy have been  compared to a SVM with a
  conventional Gaussian  kernel on the  original data and on  the data
  projected on  the first  principal axis  of the  considered classes,
  called the PCA-Mahalanobis kernel.   The main difference between the
  HDDA-MK and  the PCA-Mahalanobis kernel is  that the PCA-Mahalanobis
  kernel discards the  noise subspace while the  HDDA-MK also exploits
  the noise subspace in order to improve the class discrimination.  As
  previously  stated,  Gaussian   kernel  and  PCA-Mahalanobis  kernel
  correspond to extreme cases of HDDA-MK: $\sigma_{\hat{p}}^2=+\infty,
  \forall              p\in\{1,\dotsc,\hat{p}_{N_c}\}$              or
  \mbox{$\sigma_{\hat{p}+1}^2=+\infty$}.}

Simulated   data    were   constructed   using    a   linear   mixture
model~\cite{974727}:
\begin{eqnarray}   
  \mathbf{x}  =   \sum_{i=1}^{N_c}\alpha_i\mathbf{s}_i  +\mathbf{b}
\end{eqnarray} 
where   $N_c$   is   the   number    of   classes,   $y=j$   such   as
$\alpha_j=\max_i{\alpha_i}$,
$b\sim\mathcal{N}(\mathbf{0},\varepsilon^2\mathbf{I})$             and
$\mathbf{s}_i$   follows  the   HDDA  model.    The  mean   values  of
$\mathbf{s}_i$ were  extracted from  a spectral library  provided with
the ENVI software used  in hyperspectral imagery~\cite{envi:soft}. The
number of spectral variables $d$ was set to $413$ and $p_c$ was set to
10 for each class. The  noise variance $\varepsilon^2$ was adjusted to
get a SNR~=~1.   Three experiments were run for a  different number of
classes,   i.e.,    $2$,   $3$   and   $4$    classes,   respectively.
Figure~\ref{fig:simulated:spectra} presents two  simulated spectra and
their linear mixture.  The number of training samples was 1000 and the
number of testing samples was  $1500$.  The experiment was repeated 50
times  for each  configuration.   The  hyperparameters were  estimated
using  the  radius  margin  bound,  for  each  classifier.   Since  no
difference  in terms  of classification  accuracies were  observed, we
only report the results of the ``1 vs all'' classifier.

\subsubsection{Estimation of $\hat{p}_c$ }
\begin{figure}
  \centering
  \includegraphics[width=0.3\textwidth]{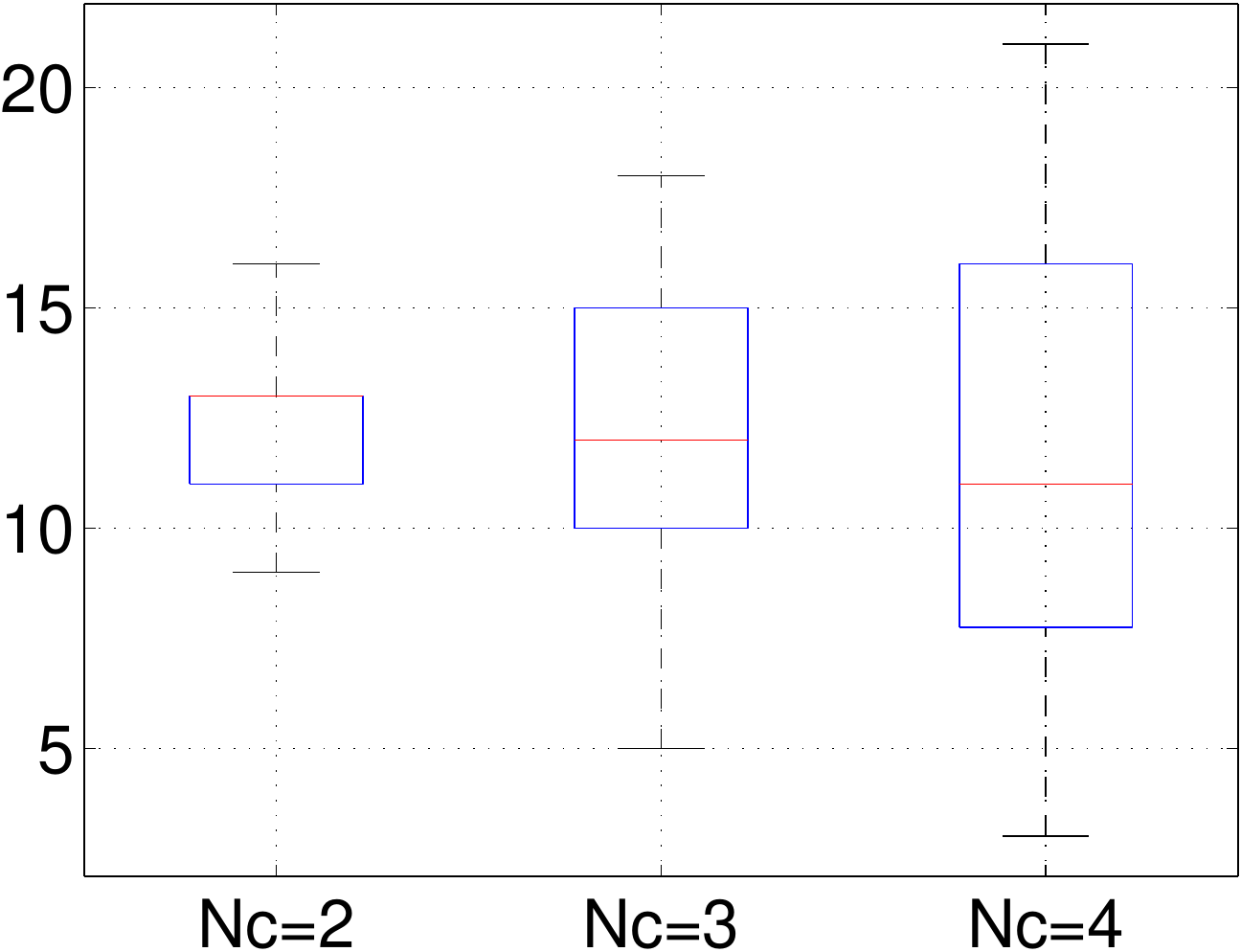}
  \caption{Boxplot  of the  estimation  of $\hat{p}_c$  for the  three
    configurations.   The   vertical   axis   represents   values   of
    $\hat{p}_c$.}
  \label{fig:p:estim}
\end{figure}

With simulated  data, it  is possible  to assess how  the size  of the
intrinsic  signal  subspace  is  estimated.   Figure~\ref{fig:p:estim}
presents the  boxplots of the  estimations.  After several  tries, the
threshold for the  scree test was fixed to  10\%. Fixing the threshold
to a  too high value would lead  to underestimate $p$ while  a too low
value would lead to drastically overestimate it.

From the figure,  the scree test overestimates the  parameter $p$, for
each configuration.  The variance of the estimation is larger when the
number  of classes  is  increased  while the  bias  of the  estimation
decreases.  However, the error in  estimating $p$ is not too important
with regards to the original size of the data ($d=413$).  Furthermore,
in    previous    work~\cite{5651956},    another    criterion    (the
BIC~\cite{BIC})  was used  to estimate  the correct  dimension of  the
subspace where the  data live.  The BIC criterion  showed poor results
when the number  of training samples for single class  $n_c$ was close
to the dimension of the  data ($d\approx n_c$).  From the experiments,
the scree test is more robust in such a situation.

\subsubsection{Classification accuracies}
\begin{figure*}
  \centering
  \begin{tabular}{ccc}
    \includegraphics[width=0.3\textwidth]{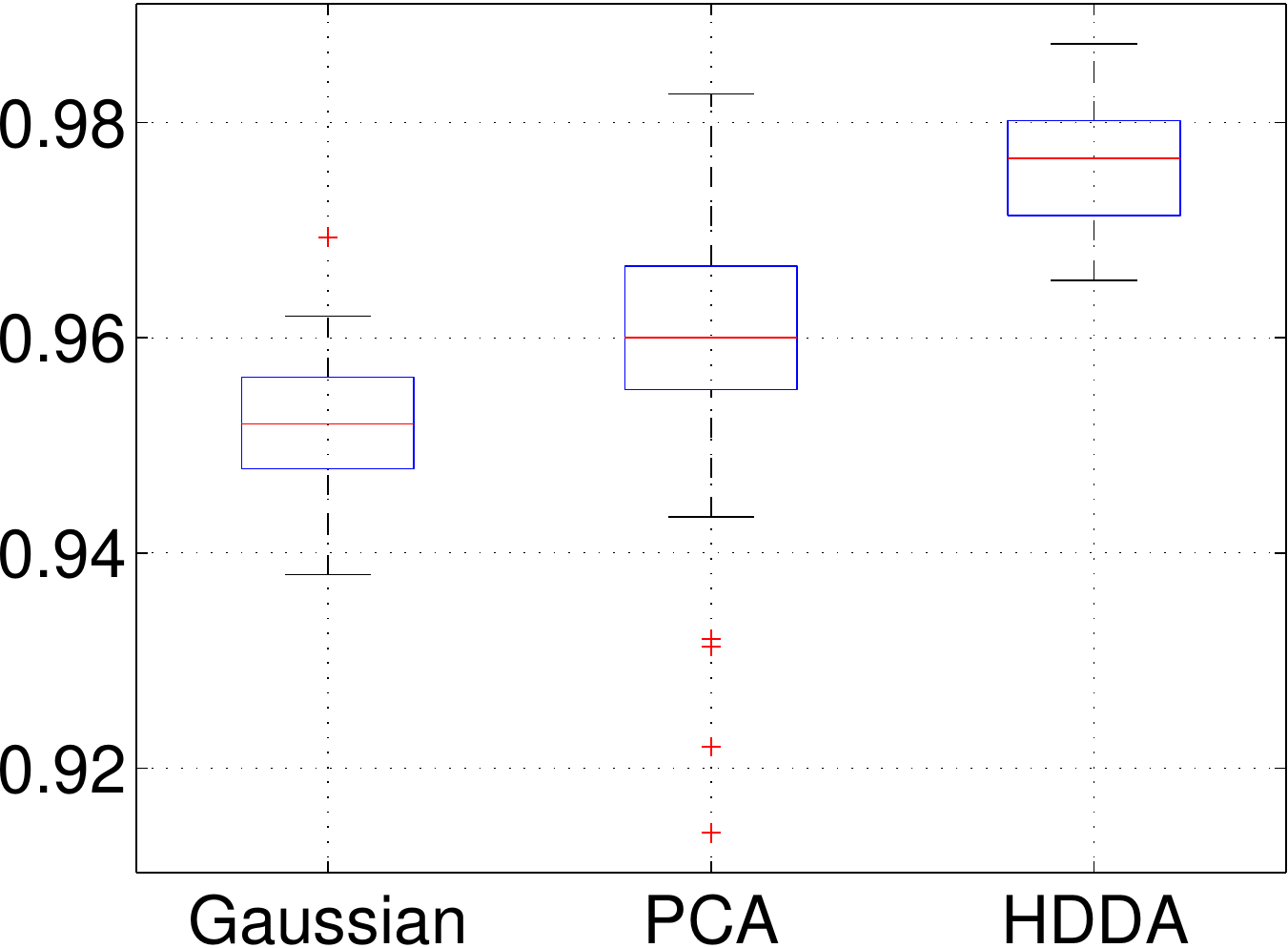} & \includegraphics[width=0.3\textwidth]{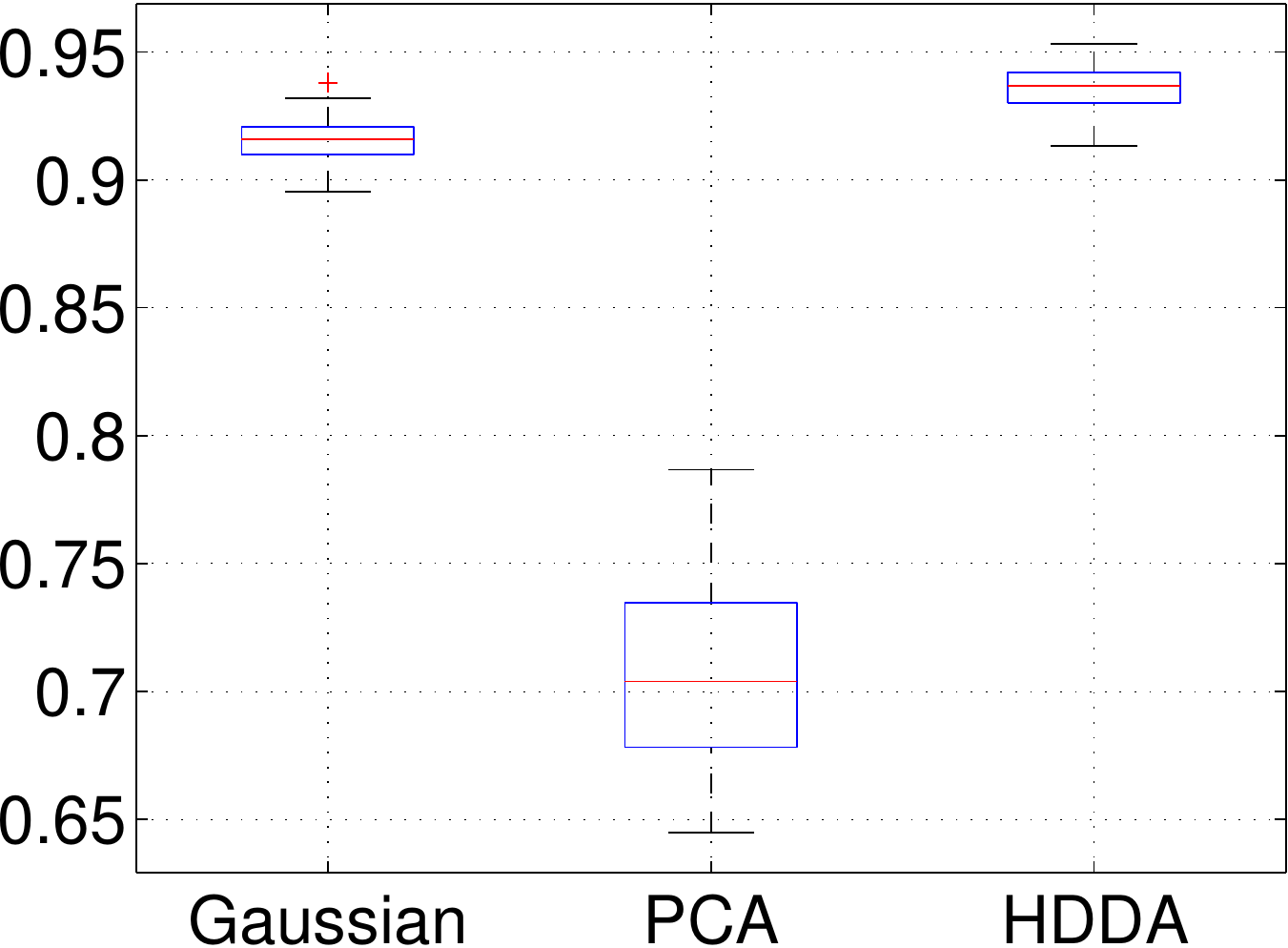} & \includegraphics[width=0.3\textwidth]{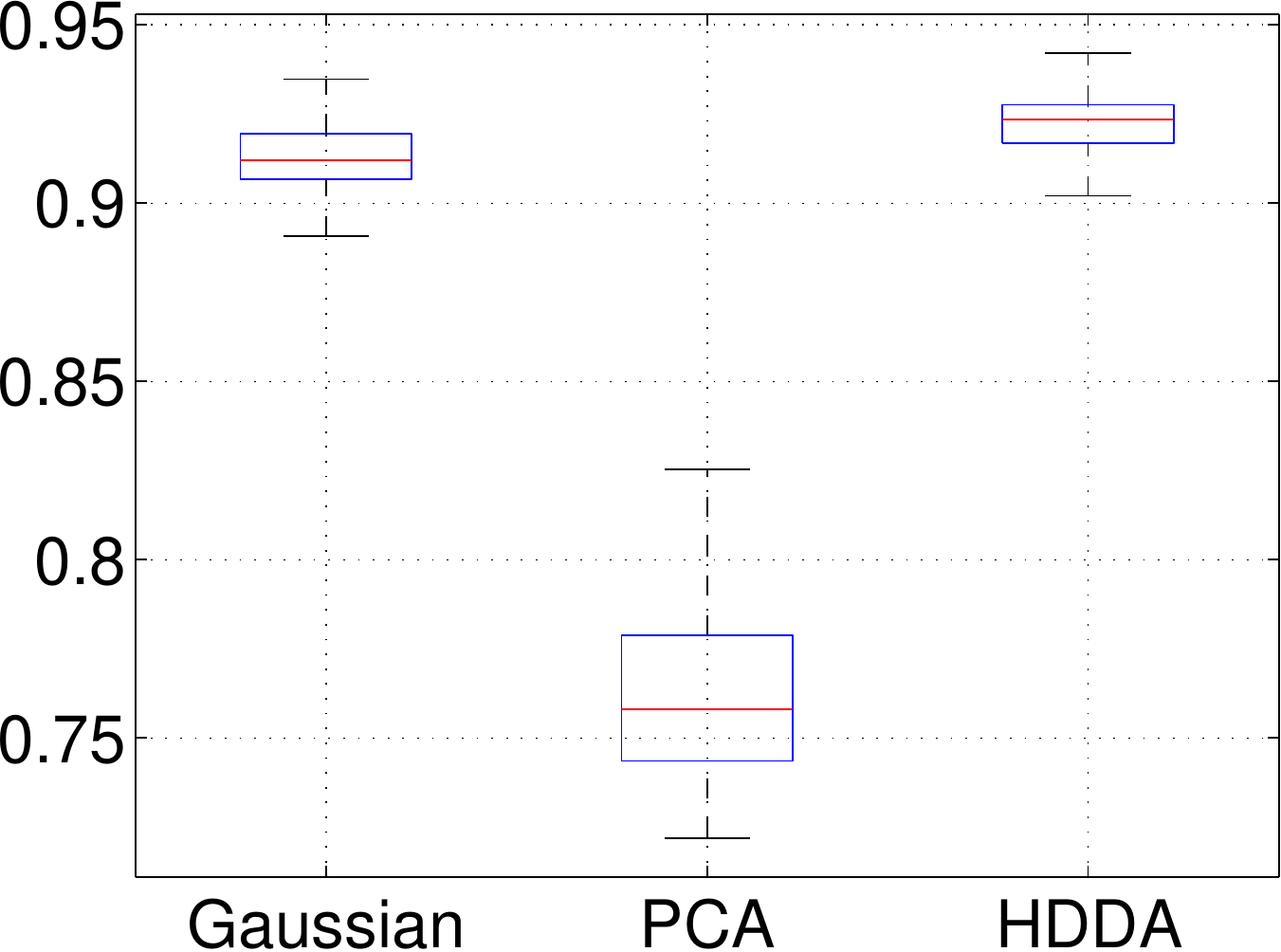}\\
    (a) & (b) & (c)
  \end{tabular}
  \caption{Boxplots  of the  classification accuracies  for  the three
    experiments.   (a)  $N_c=2$, (b)  $N_c=3$  and  (c) $N_c=4$.   The
    vertical axis represent the overall classification accuracies.}
  \label{fig:simu:oa}
\end{figure*}

The   percentages   of   correct  classification   are   reported   in
Figure~\ref{fig:simu:oa}.   For the  three  experiments, the  proposed
kernel  leads to  the best  results in  terms of  accuracies. Although
$\hat{p}_c$ was overestimated, it did not penalize the performances of
the algorithm in terms of  classification accuracies. For $N_c=2$, the
second best  result is provided  by the PCA-Mahalanobis  kernel, while
for $N_c=3$ or 4 it is provided  by the Gaussian kernel applied on the
original data.  For  instance, for $N_c=4$, the mean  value of correct
classification is  92.2\% for the HDDA-Mahalanobis  kernel, 91.3\% for
the   conventional   Gaussian  kernel   and   only   76.3\%  for   the
PCA-Mahalanobis.   The  results   confirm   the  poor   generalization
capability  of   the  Mahalanobis   kernel  when  dealing   with  high
dimensional spaces.  Although the conventional Gaussian kernel is less
sensitive  to the  problem, the  proposed kernel  gives a  significant
improvement of the classification accuracy.

\subsection{Classification of Madelon data}
Madelon  data set  is a  simulated data  set used  for the  \emph{NIPS
  Feature                                                    Selection
  Challenge}\footnote{\url{http://www.nipsfsc.ecs.soton.ac.uk/}}.   It
has 5 useful features, 15 redundant features and 480 random probes for
a total of  500 features ($d=500$). It is composed  of two classes and
the  number of  training samples  is 2000  and the  number of  testing
sample is  $600$. The  threshold for  the scree test  has been  set to
20\%. The proposed kernel has been  compared to the same kernels as in
the previous section.

The  classification results  in terms  of accuracies  are reported  in
Table~\ref{table:madelon}.  The conventional  Gaussian kernel performs
badly  on that  data set  with a  global precision  of 69,7\%  (random
classifier would achieve 50\%). The  best accuracy is obtained for the
proposed  kernel  with  an  average  accuracy  of  83.9\%.   The  size
$\hat{p}_c$  of   the  signal  subspace   for  the  two   classes  was
$\hat{p}_1=4$ and  $\hat{p}_2=4$, respectively.  The  results obtained
with the PCA-Mahalanobis kernel are worse than those obtained with the
HDDA-Mahalanobis kernel. Thus, it confirms  the pertinence to use both
the  signal  and  the  noise  subspace from  the  HDDA  model  in  the
classification.

\begin{table}
  \caption{Madelon data: Percentage of samples correctly classified.}
  \label{table:madelon}
  \centerline{\begin{tabular}{|l|c|c|c|}
      \hline
      & Class 1 & Class 2& Mean\\
      \hline
      \hline
      Gaussian& 69.7 & 69.7& 69.7\\
      \hline
      PCA-Mahalanobis & 83.3 & 81.8&82.5\\
      \hline
      HDDA-Mahalanobis &\bf 84.1 &\bf 83.8 &\bf 83.9\\
      \hline
    \end{tabular}}
\end{table}

\subsection{Classification of Arcene data set}
\textcolor{black}{Arcene data set  is data set used  for the \emph{NIPS
    Feature Selection  Challenge}.  It  has 7000 real  variables, 3000
  random  probes for  a total  of  10000 features  ($d=10000$). It  is
  composed of two classes. The number of available training samples is
  100 and the number of test samples is 100.  Therefore, the number of
  training  samples is  very small  in comparison  with the  number of
  variables.  About 50\% of the data are non zero.}

\textcolor{black}{  The  classification   accuracies  are  reported  in
  Table~\ref{tab:arcene}.  The  selected threshold for the  scree test
  has been set to 0.5\%. The  results obtained for other values of the
  threshold  are also  reported for  comparison.  The  Gaussian kernel
  achieves a global accuracy of 80\%.  It performs quite well in terms
  of classification accuracies  related to the dimension  of the data.
  The  PCA-Mahalanobis kernel  performs worst  whatever the  threshold
  value.  For  the HDDA-Mahalanobis kernel,  for the highest  value of
  $s$ the results are equal to  those obtain with the Gaussian kernel.
  Then a slight  increased of the accuracy is  observed for $s=0.005$.
  When the  size of the signal  subspace is too large  ($s=0.0001$ and
  ($\hat{p}_1,\hat{p}_2)=(22,23)$), too  many hyperparameters  have to
  be estimated and thus the classification accuracy becomes low.}

\begin{table}\small
  \caption{Arcene data: $s$ is the threshold value in the scree test, $(\hat{p}_1,\hat{p}_2)$ correspond to the estimated size of  the signal subspace for each class.  The number is the percentage of samples correctly classified.}
  \label{tab:arcene}
  \begin{tabular}{|c|c|cccccc|}
    \hline
    &$s$ & 0.2 & 0.1 & 0.05 & 0.01 & 0.005 & 0.001\\
    \hline
    &$(\hat{p}_1,\hat{p}_2)$ & (2,2) & (2,3) & (3,4) & (8,7) & (11,10) & (22,23)\\
    \hline
    \multirow{2}{*}{PCA-Mahalanobis} & Class 1 & 69.0 & 69.0 & 70.0 & 75.0 & 78.0 & 51.0\\
    & Class 2 & 72.0 & 72.0 & 62.0 & 67.0 & 72.0 & 75.0\\
    \hline
    \multirow{2}{*}{HDDA-Mahalanobis} & Class 1 &  80.0 & 80.0 & 80.0 & 80.0 &\bf 83.0 & 51.0\\
    & Class 2 &  80.0 & 80.0 &\bf 81.0 & 80.0 &\bf 81.0 & 74.0\\
    \hline
  \end{tabular}
\end{table}
\subsection{Classification of real hyperspectral data}
The  data set considered  in this  experiment is  the \emph{University
  Area} of Pavia, Italy, acquired with the ROSIS-03 sensor.  The image
has  103 spectral  variables, i.e.,  each  pixel is  represented by  a
vector  with 103  features (d=103)~\cite{Chen20082731}.   Nine classes
have been defined  by photo-interpretation as seen in  first column of
Table~\ref{tab:res}.  Here,  the threshold was set  to 0.01\%, because
of a very  high value of the first principal  component (mainly due to
the albedo).

Classification  results  are  reported  in  Table~\ref{tab:res}.   The
proposed kernel leads to an  increase of the accuracy, compared to the
conventional  Gaussian  kernel.   However,  for  this  data  set,  the
PCA-Mahalanobis  and HDDA-Mahalanobis kernel  perform equally  well in
terms  of  accuracies,  except   for  the  classes  \emph{meadow}  and
\emph{bare soil}.

\begin{figure}
  \centering
  \includegraphics[width=0.3\textwidth]{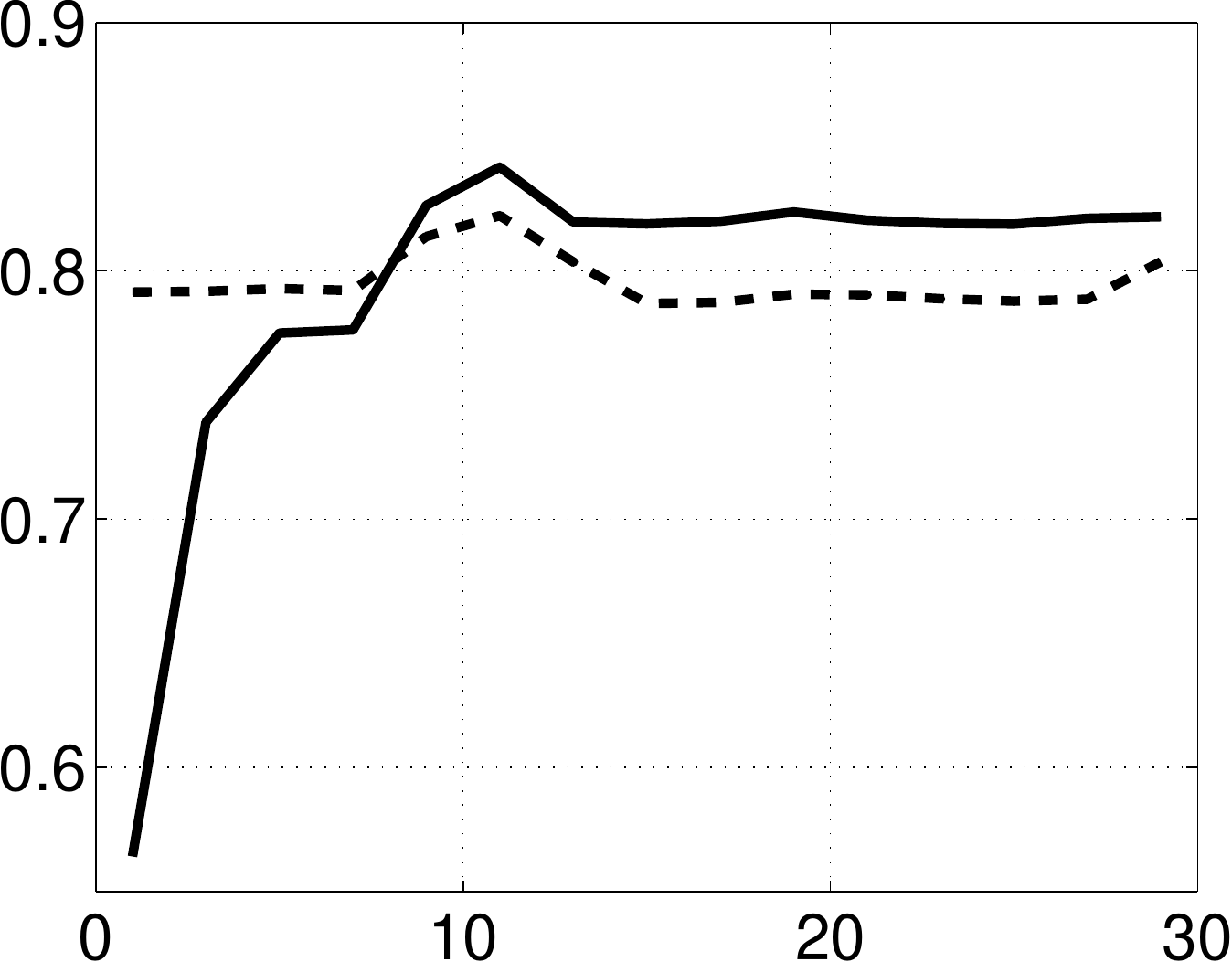}
  \caption{Classification accuracies as a function of the dimension of
    the signal  subspace for the class  \emph{meadow}.  The horizontal
    axis  represents  the  value  $\hat{p}_c$ and  the  vertical  axis
    represents   the   classification   accuracies   for   the   class
    \emph{meadow}.  Solid  line  corresponds  to  the  PCA-Mahalanobis
    kernel  and the  dashed line  corresponds to  the HDDA-Mahalanobis
    kernel.  The selected value with the scree test is $10$.}
  \label{fig:p}
\end{figure}

To  assess   the  influence  of  $\hat{p}_c$   on  the  classification
accuracies, the class \emph{meadow} has been classified for a range of
values   of   $\hat{p}_c$.    The   results  are   reported   in   the
Figure~\ref{fig:p}   for  the  PCA-Mahalanobis   and  HDDA-Mahalanobis
kernels.   From the  figure, the  optimal $\hat{p}_c$  is 11  which is
close to the value selected  with the scree test ($\hat{p}_c=10$). The
cumulative  variance  is   99.72\%  for  $\hat{p}_c=10$,  99.75\%  for
$\hat{p}_c=11$ and 99.77\% for $\hat{p}_c=12$.  The proposed kernel is
slightly  influenced by  the  choice of  $\hat{p}_c$  for that  class,
compared to PCA-Mahalanobis kernel. In particular, the proposed kernel
is robust  if $\hat{p}_c$ is underestimated.   However, if $\hat{p}_c$
is  heavily  overestimated,  too  many  parameters would  need  to  be
estimated and it could degrade the training process.

\begin{table*}
  \small
  \centering
  \caption{Classification accuracies for the different kernels in percentage of correctly classified samples. In the first column, the numbers in brackets represent the total number of training and testing samples for each class, respectively.}
  \label{tab:res}
  \begin{tabular}{|l|c||c|c|c|}
    \hline
    \multicolumn{1}{|c}{}& $\hat{p}_c$& Gaussian & PCA-Mahalanobis & HDDA-Mahalanobis\\
    \hline
    \hline
    Asphalt (548, 6631) & 12 &94.8 &\bf 95.8&\bf 95.8 \\
    \hline
    Meadow (540, 18649) & 10 &79.4 &\bf 83.6&82.1 \\
    \hline
    Gravel (392, 2099)  & 9 & 97.2 &\bf 97.5&97.2 \\
    \hline
    Tree (524, 3064) & 14 &94.3 & \bf 98.2&\bf 98.2 \\
    \hline
    Metal Sheet (265, 1345) & 7 &99.8 &\bf 99.9&\bf 99.9\\
    \hline
    Bare Soil (532, 5029) & 9 &87.8 &85.9&\bf 88.4 \\
    \hline
    Bitumen (375, 1330) & 21 &98.8 &98.7&\bf 99.0 \\
    \hline
    Brick (514, 3682) & 12 &96.7 &\bf 97.2&\bf 97.2\\
    \hline
    Shadow (231, 947) & 14 &\bf 99.9&\bf 99.9&\bf 99.9\\
    \hline
    \hline
    \multicolumn{2}{|c||}{Average class accuracy}&94.3 &95.2&\bf 95.3\\
    \hline
  \end{tabular}
\end{table*}

\subsection{Analysis of the processing time}
\textcolor{black}{To  assess  the computational  load  of the  proposed
  method, the processing time was computed for four data sets. The two
  data sets from the \emph{NIPS Feature Selection Challenge} were used
  and  well as  the  two first  classes  (Asphalt and  Method) of  the
  hyperspectral    data   set.   The    results   are    reported   in
  Table~\ref{tab:tp}.  The program runs  under Matlab  on a  two cores
  2.67GHz laptop.}

\textcolor{black}{For the  Arcene data  set, the Gaussian  kernel shows
  the  lowest  computational  time.   The  computation  of  the  first
  eigenvalues/eigenvectors  is  demanding since  the  data have  10000
  features.  It  requires about 14  seconds.  The optimization  of the
  hyperparameters is fast  since a few number of  training samples are
  available. For  the Madelon data  sets, the computation of  the firs
  eigenvalues/eigenvectors  is  fast,  about  1.3  seconds  while  the
  optimization of  the kernel hyperparameters is  more demanding.  For
  that data  set, the HDDA-Mahalanobis is the  slowest.  Regarding the
  University data  set, for the  first two classes, the  estimation of
  the first  eigenvalues/eigenvectors is very fast,  about 0.3 second.
  However,  the  estimation  of  the kernel  hyperparameters  is  more
  demanding  than with  the  two  others data  sets.  For the  Asphalt
  problem, the HDDA-Mahalanobis performs  the fastest, while it is the
  slowest for the Meadow problem.}

\textcolor{black}{From the above results,  it is difficult to point out
  a clear  winner in terms  of processing time.  From  a computational
  viewpoint, the optimization of the Gaussian kernel is less demanding
  than   PCA-  or   HDDA-Mahalanobis  kernels.   However,   since  the
  optimization problem is a gradient descent on a non-convex function,
  a local  minimum might be found  sooner with one method  and make it
  faster  than  the  others.  Nevertheless,  the  HDDA-Mahalanobis  is
  usually more  demanding in  terms of processing  time. Computational
  complexity  of  PCA-  and  HDDA-Mahalanobis  can be  assumed  to  be
  comparable in terms of processing time.}

\begin{table}
  \centering
  \caption{Processing time in seconds of the competive methods for four different data sets.}
  \label{tab:tp}
  \begin{tabular}{|c||c|c|c|c|}
    \hline
    & Arcene & Madelon & Asphalt & Meadow\\
    \hline
    Gaussian &\bf 2.4 & 31.9 & 113 .0 & \bf 87.4\\
    PCA-Mahalanobis & 14.2 &\bf 24.0 & 193.0 & 185.7\\
    HDDA-Mahalanobis & 14.5 & 80.8 &\bf 94.2 & 195.3\\
    \hline
    \end{tabular}    
\end{table}
\section{Conclusions}
\textcolor{black}{In  this  paper,  a  novel  kernel  adapted  to  high
  dimensional  data has been  proposed.  The  parsimonious Mahalanobis
  kernel is based on the  \emph{emptiness} property of HD spaces.  For
  each class, the original input space is split into a signal subspace
  and a noise  subspace.  Using this assumption, the  inversion of the
  covariance  matrix  in  the  Mahalanobis kernel  can  be  accurately
  computed.  The proposed kernel was tested in a SVM framework for the
  purpose of  classification.  Experimental results on  four data sets
  demonstrate the potential of the proposed kernel.  In each case, the
  classification  accuracy  increased  compared  to  the  conventional
  Gaussian  kernel and  for  three cases  the  proposed kernel  showed
  superior  results to  simply map  the data  on the  first  PCA axes.
  Consequently,  for HD  data  the HDDA-Mahalanobis  kernel should  be
  prefered.}

Regarding the computational load,  the HDDA-Mahalanobis kernel is more
demanding during  the training process than the  Gaussian kernel since
more  hyperparameters  have  to   be  estimated.   Besides  that,  the
HDDA-Mahalanobis kernel  is efficient when the dimension  of the input
space is high.  Hence,  for moderate or small dimensions, conventional
kernels should be preferred.

In  the  article,  only  classification is  investigated.   But  other
processing     could    also     have    been     considered,    e.g.,
regression~\cite{abe:regression}. The critical point  is to be able to
tune  the  hyperparameter  efficiently,   which  it  is  feasible  for
regression. However, the actual optimization of the hyperparameters is
demanding  in terms  of  computations  and it  is  sensitive to  local
minima. Therefore, a different strategy must be studied.

Perspectives of this work concern the development of new kernels using
the  HDDA-model.  For instance,  it is  possible to  define a  new dot
product for the conventional polynomial kernel:
\begin{eqnarray}k(\mathbf{x},\mathbf{z}|c)                            =
  \big(\mathbf{x}^t\boldsymbol{\Sigma}_c^{-1}\mathbf{z}+1\big)^r.\end{eqnarray}
Furthermore,  a mixture  of  kernels  using the  HDDA  could be  used.
From~(\ref{eq:mixture:prod}),  the   HDDA-Mahalanobis  kernel  can  be
considered as  a product of several  kernels.  In the  future, we will
investigate the summation of kernels. 

Finally, a free-parameter alternative to the scree test for estimation
of the intrinsic dimension must be addressed. For instance, an maximum
likelihood    estimator     for    HDDA    exits     and    must    be
investigated~\cite{Bouveyron20111706}.

\section*{Acknowledgment}
The authors would like to thank  the IAPR - TC7 for providing the data
and Prof. Paolo Gamba and  Prof. Fabio Dell'Acqua of the University of
Pavia,  Italy,  for  providing   reference  data. 

\bibliography{IEEEabrv,bib}
\end{document}